\documentclass{article}
\usepackage{a4wide}
\usepackage{amssymb}
\usepackage{amsfonts}
\usepackage{amsmath}
\usepackage{amsxtra}
\usepackage{latexsym}
\usepackage{mathabx}
\usepackage{mathdots}
\usepackage{graphicx}
\usepackage{hyperref}
\usepackage{authblk}
\usepackage{tikz}
\usepackage{xcolor}
\usepackage{xurl}
\usepackage{fancyvrb}
\usepackage{lipsum}
\usepackage{comment}
\usepackage[title]{appendix}

\DeclareMathOperator\inv{inv}

\DeclareMathOperator{\HW}{HW}

\DeclareMathOperator{\ord}{ord}
\DeclareMathOperator{\dlog}{dlog}

\begin{document}

\newtheorem{theo}{Theorem}[section]
\newtheorem{theo*}{Theorem}
\newtheorem{lemma}[theo]{Lemma}
\newtheorem{definition}[theo]{Definition}
\newtheorem{notation}[theo]{Notation}
\newtheorem{obs}[theo]{Observation}
\newtheorem{remark}[theo]{Remark}
\newtheorem{cor}[theo]{Corollary}
\newtheorem{cor*}{Corollary}
\newtheorem{prop}[theo]{Proposition}
\newtheorem{conj}[theo]{Conjecture}
\newtheorem{example}[theo]{Example}

\newcommand{\N}{{\mathbb{N}}}
\newcommand{\Z}{{\mathbb{Z}}}  
\newcommand{\Q}{{\mathbb{Q}}}
\newcommand{\R}{{\mathbb{R}}}
\newcommand{\C}{{\mathbb{C}}}
\newcommand{\M}{{\mathbb{M}}}
\newcommand{\F}{{\mathbb{F}}}
\newcommand{\V}{{\mathbb{V}}}
\newcommand{\LL}{{\mathbb{L}}}
\newcommand{\A}{{\mathbb{A}}}
\newcommand{\HH}{{\mathbb{H}}}

\newcommand{\g}[1]{{\mathfrak {#1}}}
\newcommand{\cc}[1]{{\cal {#1}}}
\newcommand{\ul}{\underline }
\newcommand{\ol}{\overline }

\newcommand{\qed}{\hspace*{\fill}$\Box$}
\newcommand{\Qed}{\hspace*{\fill}$\Box \Box$}

\title{On the representation of number-theoretic functions by arithmetic terms}

\author{Mihai Prunescu \footnote{Research Center for Logic, Optimization and Security (LOS), Faculty of Mathematics and Computer Science, University of Bucharest, Academiei 14, Bucharest (RO-010014), Romania; and Simion Stoilow Institute of Mathematics of the Romanian Academy, Research unit 5, P. O. Box 1-764, Bucharest (RO-014700), Romania. E-mail: {\tt mihai.prunescu@imar.ro}, {\tt mihai.prunescu@gmail.com}.}, Lorenzo Sauras-Altuzarra \footnote{Simion Stoilow Institute of Mathematics of the Romanian Academy, Research unit 5, P. O. Box 1-764, Bucharest (RO-014700), Romania. E-mail: {\tt lorenzo@logic.at}. Partially supported by Bitdefender (Postdoctoral Fellowship).}}

\date{}

\maketitle

\begin{abstract}

We present closed forms for several functions that are fundamental in number theory and we explain the method used to obtain them. Concretely, we find formulas for the $ p $-adic valuation, the number-of-divisors function, the sum-of-divisors function, Euler's totient function, the modular inverse, the integer part of the root, the integer part of the logarithm, the multiplicative order and the discrete logarithm. Although these are very complicated, they only involve elementary operations, and to our knowledge no other closed form of this kind is known for the aforementioned functions.

\vspace{1mm}


\vspace{1mm}


\end{abstract} 

\section{Introduction}\label{SectionIntroduction}

For any non-empty set $ X $ and any non-empty set $ F $ of finitary operations on $ X $, we define the \textbf{inductive closure} of $ X $ with respect to $ F $ as the minimum set $ C $ for which $ X \subseteq C $ and such that, if $ r $ is a positive integer, $ f $ is an $ r $-ary operation in $ F $ and $ \vec{c} \in C^r $, then $ f ( \vec{c} ) \in C $ (cf.\ Enderton \cite[Section 1.4]{Enderton}).

We denote the set of non-negative integers by $ \mathbb{N} $.

The binary operation on $ \mathbb{N} $ that is given by $ \max ( x - y , 0 ) $ is called \textbf{truncated subtraction} and denoted by $ \dotdiv $ (see Vereschchagin \& Shen \cite[p.\ 141]{VereschchaginShen}).

For any positive integer $ r $, we define an $ r $-variate \textbf{arithmetic term} in variables $ n_1 $, $ \dots $, $ n_r $ as an element of the inductive closure of $ \mathbb{N} \cup \{ n_1 , \dots , n_r \} $ with respect to the binary operations on $ \mathbb{N} $ given by $$ x + y , x \dotdiv y , \left\lfloor x / y \right\rfloor , x^y , $$ (cf.\ Prunescu \& Sauras-Altuzarra \cite{PrunescuSaurasAltuzarra}). We follow the conventions $ 0^0 = 1 $ (see Mendelson \cite[Proposition 3.16]{Mendelson}) and $ \left\lfloor x / 0 \right\rfloor = 0 $ (see Mazzanti \cite[Section 2.1]{Mazzanti}).

Note that the total number of operations occurring in an arithmetic term is \textit{fixed} (i.e.\ it does not depend on the arguments). The expressions satisfying this condition are usually known as \textit{closed forms} (cf.\ Borwein \& Crandall \cite{BorweinCrandall}).

For example, $ 2^{n + 1} - 1 $ and $ \sum_{k = 0}^n ( 2^k ) $ are expressions that represent the same integer sequence. However, only the first one is considered a closed form, because the total number of operations of the latter depends on the argument $ n $.

A celebrated kind of closed form is the so-called \textbf{hypergeometric closed form}: a linear combination, with respect to a field $ K $, of expressions $ f ( n ) $ such that $ f ( n + 1 ) / f ( n ) $ is a rational function on $ K $ (cf.\ Petkovšek et al.\ \cite[Definition 8.1.1]{PetkovsekEtAl} and Sauras-Altuzarra \cite[Definition 1.4.13]{SaurasAltuzarra}).

A \textbf{Kalmar function} is a computable and finitary operation on $ \mathbb{N} $ whose deterministic computation time is upper-bounded by some integer sequence of the form $$ 2^{2^{\iddots^{2^n}}} $$ (cf.\ Marchenkov \cite[Introduction]{Marchenkov}, Prunescu \& Sauras-Altuzarra \cite{PrunescuSaurasAltuzarra} and Oitavem \cite[Introduction]{Oitavem}).

In most of the mathematical contexts, the concrete integer sequences of non-negative terms that one encounters are Kalmar functions. Amazingly, Mazzanti \cite{Mazzanti}, and later Marchenkov \cite{Marchenkov}, achieved to prove that every Kalmar function admits an arithmetic-term representation (of the same number of arguments).

As the proofs displayed by Mazzanti \cite{Mazzanti} and Marchenkov \cite{Marchenkov} involve some special kinds of exponential Diophantine representations, the potential applications in number theory remained unobserved. In Section \ref{SectionRepresentationMethod}, we reformulate their arithmetic-term representation method in order to make its applicability explicit.

Given two coprime integers $ n \geq 2 $ and $ m \in \{ 1 , \dots , n - 1 \} $, we define the number $ \inv ( m , n ) $ as the \textbf{modular inverse} of $ m $ modulo $ n $ (i.e. as the only number $ x \in \{ 1 , \dots , n - 1 \} $ such that $ m x \equiv 1 \pmod{n} $, cf.\ Rosen \cite[Section 4.2]{Rosen}). And we define the number $ \ord ( m , n ) $ as the \textbf{multiplicative order} of $ m $ modulo $ n $ (i.e. as the minimum positive integer $ r $ such that $ m^r \equiv 1 \pmod{n} $, cf.\ Rosen \cite[Section 9.1]{Rosen}).

Recall that \textbf{Euler's totient function} $ \varphi $ counts the positive integers that do not exceed its argument $ n $ while being coprime with $ n $ (see K\v{r}\'{i}\v{z}ek et al.\ \cite[Table 13.3]{KrizekEtAl}).

Given two coprime integers $ n \geq 2 $ and $ g \in \{ 1 , \dots , n - 1 \} $, the number $ g $ is said to be a \textbf{primitive root} modulo $ n $ if, and only if, $ \ord ( g , n ) = \varphi ( n ) $ (cf.\ Rosen \cite[Section 9.1]{Rosen}). If this is the case, then, given also another integer $ m \in \{ 1 , \dots , n - 1 \} $ that is coprime with $ n $, the \textbf{discrete logarithm} of $ m $ to the base $ g $ modulo $ n $, which we denote by $ \dlog ( m , g , n ) $, is the only number $ d \in \{ 1 , \ldots , \varphi ( n ) \} $ such that $ g^d \equiv m \pmod{n} $ (cf.\ Crandall \& Pomerance \cite[Subsection 6.4.1]{CrandallPomerance} and Rosen \cite[Section 9.4]{Rosen}).

After Section \ref{SectionUsefulTerms} displaying arithmetic terms that we repeatedly use and Section \ref{SectionRepresentationMethod} explaining the aforementioned method for arithmetic-term construction, some examples are presented. Namely, we apply this method to the number-of-divisors function $ \tau $ in Section \ref{SectionTau}, the sum-of-divisors function $ \sigma $ in Section \ref{SectionSigma}, Euler's totient function $ \varphi $ in Section \ref{SectionPhi}, the modular inverse $ \inv $ in Section \ref{SectionInv}, the integer part of the root in Section \ref{SectionRoot}, the integer part of the logarithm in Section \ref{SectionLog}, the multiplicative order $ \ord $ in Section \ref{SectionOrd} and the discrete logarithm $ \dlog $ in Section \ref{SectionDiscLog}. The correctness of many of these calculations can be experimentally verified with the Maple codes from Appendix \ref{MapleCodes}.

The above functions are usually computed by algorithms that inspect the numbers from $ 1 $ to $ n $. Also, if the prime factor decomposition of $ n $ is $ n = p_1^{k_1} \dots p_r^{k_r} $, where $ p_1 , \dots , p_r $ are pairwise distinct primes and $ k_1 , \dots , k_r $ are positive integers, then \begin{eqnarray} \tau(n) & = & (k_1 + 1) \dots (k_r + 1) , \\ \sigma(n) & = & \frac{p_1^{k_1+1} - 1}{p_1 - 1} \dots \frac{p_r^{k_r+1} - 1}{p_r - 1} , \\ \varphi(n) & = & n \left( 1 - \dfrac{1}{p_1} \right) \dots \left( 1 - \dfrac{1}{p_r} \right) \label{EqCharPhi} \end{eqnarray} (see Hardy \& Wright \cite[Theorem 275, Theorem 273 and Theorem 62]{HardyWright}). Thus, at a first sight, it is difficult to imagine that one can compute these functions by applying a fixed number of arithmetic operations in a given order to the number $ n $ alone (in particular, without a prior knowledge of its prime number decomposition), but arithmetic terms do precisely this task.

The arithmetic terms obtained in the present article are however very long and complicated, and therefore not of practical use. We believe that in the future shorter arithmetic terms will be found but, for now, we just emphasize the existence of a method to obtain them. In addition, it will be a very interesting challenge to prove (if possible, in a constructive way) that a given arithmetic term is the shortest one (among all those which compute the same Kalmar function).

\section{Useful arithmetic terms}\label{SectionUsefulTerms}

The first fundamental operation with which we will enlarge our set of admissible arithmetic terms is the product, because we have that $$ n m = \left\lfloor \dfrac{2^{n + m + 4}}{\left\lfloor \left\lfloor \dfrac{2^{n + m + 4}}{n + 1} \right\rfloor / ( m + 1 ) \right\rfloor} \right\rfloor \dotdiv ( n + m + 1 ) $$ (see Marchenkov \cite[Section 2]{Marchenkov}).

It is immediate that for every two non-negative integers $ m $ and $ n $ we have that $$ n \bmod m = n \dotdiv ( m \lfloor n / m \rfloor ) , $$ so the arithmetic term $ n \bmod m $ will be also used in representations. Observe that $ n \bmod 0 = n $ (recall that $ \lfloor n / 0 \rfloor = 0 $) and $ n \bmod 1 = 0 $.

Another admissible operation that we will need in some proofs is the maximum: $$ \max ( m , n ) = \lfloor ( m + n + ( m \dotdiv n ) + ( n \dotdiv m ) ) / 2 \rfloor . $$

Another very useful identity is: \begin{equation}\label{EqMarchenkov} n^m = 2^{( n m + n + 1 ) m} \bmod ( 2^{n m + n + 1} \dotdiv n ) \end{equation} (see Marchenkov \cite[Section 2]{Marchenkov}).

The \textbf{$ \boldsymbol{p} $-adic valuation} of $ n $ (when $ n $ is positive), which is denoted by $ \nu_p ( n ) $, is the exponent of $ p $ in the prime number decomposition of $ n $ (see the Encyclopedia of Mathematics \cite{Encyclopedia}).

\begin{theo}\label{ThmPadicVal} The function $ \nu_p ( n ) $ (for integer arguments $ n \geq 1 $ and $ p $ prime) can be represented by the arithmetic term $$ \left\lfloor \dfrac{\gcd(n, p^n)^{n+1} \bmod (p^{n+1} \dotdiv 1)^2}{p^{n+1} \dotdiv 1} \right\rfloor . $$ \end{theo}

{\bf Proof} Let $ x = \nu_p ( n ) $.

It is clear that $ x < n < p^n $, so $ x < p^{n + 1} \dotdiv 1 $ and consequently $ 1 + x ( p^{n + 1} \dotdiv 1 ) < ( p^{n + 1} \dotdiv 1 )^2 $.

In addition, $$ \gcd ( n , p^n )^{n + 1} = {\left ( p^x \right )}^{n+1} = {\left ( p^{n+1} \right )}^{x} =  {\left ( p^{n+1} - 1 + 1 \right )}^{x} = \sum_{k=0}^x \left( \binom{x}{k}(p^{n+1} \dotdiv 1)^k \right) = $$ $$ 1 + x (p^{n+1} \dotdiv 1) + ( p^{n+1} \dotdiv 1 )^2 \sum_{k = 2}^x \left( \binom{x}{k}( p^{n+1} \dotdiv 1 )^{k-2} \right) . $$

Thus $ \gcd(n, p^n)^{n+1} \bmod (p^{n+1} \dotdiv 1)^2 = 1 + x (p^{n+1} \dotdiv 1) $, from which the statement immediately follows. \qed

Theorem \ref{ThmEfficientPadicVal} will show an arithmetic term that computes the $ p $-adic valuation in a faster way.

Mazzanti \cite[Lemma 4.2]{Mazzanti} gave an arithmetic-term representation for the \textbf{Hamming weight} of $ n $, that is, the number of digits that are equal to one in the binary representation of $ n $ (see \href{https://oeis.org/A000120}{\texttt{OEIS A000120}}). Like later Marchenkov \cite[Section 3]{Marchenkov}, he denotes this function with $\sigma(n)$, but this notation clashes with the usual notation for the sum-of-divisors function, so in this work we replace it with $\HW(n)$.

In order to construct an arithmetic term representing $\HW(n)$, Mazzanti \cite[Lemma 3.3]{Mazzanti} proved first that $$\gcd(m,n) = \left\lfloor \dfrac{\left( 2^{m^2n (n+1)} \dotdiv 2^{m^2n} \right) \left( 2^{m^2n^2} \dotdiv 1 \right)}{\left( 2^{m^2n} \dotdiv 1 \right) \left( 2^{mn^2} \dotdiv 1 \right) 2^{m^2n^2}} \right\rfloor \bmod 2^{mn} $$ and obtained the instance of Theorem \ref{ThmPadicVal} in which $ p = 2 $. Then, by using an arithmetic-term representation of the central binomial coefficients, $$\binom{2n}{n} = \left\lfloor \dfrac{\left( 1 + 2^{2n} \right)^{2n}}{2^{2n^2}} \right\rfloor \bmod 2^{2n} , $$ and applying \textbf{Kummer's theorem} (cf.\ Matiyasevich \cite[Appendix]{Matiyasevich}), which asserts that $ \HW ( n ) $ is equal to the dyadic valuation of $ \binom{2n}{n} $, he concluded that $\HW(n)$ has a representation as an arithmetic term (which is displayed in Appendix \ref{MapleCodes}).

For every three integers $ q > 1 $, $ r \geq 0 $ and $ t \geq 0 $, there are further useful arithmetic terms representing the so-called \textbf{generalized geometric progression} of the $ r $-th kind (cf.\ Matiyasevich \cite[Appendix]{Matiyasevich}): $$ G_r (q, t) = \sum_{k = 0}^t ( k^r q^k ) . $$ Indeed, while it is well-known that \begin{equation}\label{EqGeomSum} G_0(q, t) = \dfrac{q^{t+1} - 1}{q - 1} , \end{equation} all the further $G_i(q, t)$ are recurrently obtained from the identity $$ G_{r + 1} ( q , t ) = \frac{\partial}{\partial q} ( G_r ( q , t + 1 ) ) - \sum_{j=0}^r \left( \binom{r+1}{j} G_j ( q , t ) \right) , $$ as shown by Matiyasevich \cite[Appendix]{Matiyasevich}.

We will make concrete use of the arithmetic terms $ G_1(q,t) $ and $ G_2(q, t) $, which equal $$ \dfrac{t q^{t+2} - (t+1)q^{t + 1} + 1}{( q - 1 )^2} , $$ $$ \dfrac{t^2 q^{t+3} - (2 t^2 + 2t - 1)q^{t+2} + (t+1)^2 q^{t+1} - q^2 - q}{( q - 1 )^3} , $$ respectively. These two arithmetic terms, together with $ G_0(q, t) $ and the Hamming weight, will be sufficient to represent any Kalmar function, as Corollary \ref{CorFourArithmeticTerms} will show.

Notice that all the subtractions that appear in the formulas of the generalized geometric progressions should be also written as truncated subtractions, but we left them like that in order to make the notation less cumbersome. From now on, we will keep writing the usual subtraction everywhere (in some cases for the same reason, and in the other ones because of actual need).

\section{The representation method}\label{SectionRepresentationMethod}

We call \textbf{algebraic sum} of arithmetic terms to any sum in which the summands are arithmetic terms or opposites of arithmetic terms.

For example, $ 2^x - \lfloor n / y \rfloor $ is an algebraic sum of arithmetic terms.

Given a positive integer $ k $, a \textbf{simple}-in-$ ( x_1 , \dots , x_k ) $ \textbf{exponential monomial} is an arithmetic term of the form \begin{equation}\label{ExprMonomial} \alpha ( \vec{n} ) {b_1 ( \vec{n} )}^{\beta_1 ( \vec{n} ) x_1} \dots {b_k ( \vec{n} )}^{\beta_k ( \vec{n} ) x_k} x_1^{\gamma_1} \dots x_k^{\gamma_k} , \end{equation} where $ \alpha ( \vec{n} ) $, $ \beta_1 ( \vec{n} ) $, $ \dots $, $ \beta_k ( \vec{n} ) $ are arithmetic terms, $ b_1 ( \vec{n} ) $, $ \dots $, $ b_k ( \vec{n} ) $ are positive arithmetic terms and $ \gamma_1 $, $ \dots $, $ \gamma_k $ are non-negative integers. And a \textbf{simple}-in-$ ( x_1 , \dots , x_k ) $ \textbf{exponential polynomial} is an algebraic sum of simple-in-$ ( x_1 , \dots , x_k ) $ exponential monomials.

The main goal of this section is to explain the technique of construction of arithmetic terms developed by Matiyasevich \cite[Section 6.3]{Matiyasevich}, and used also by Marchenkov \cite{Marchenkov} and Mazzanti \cite{Mazzanti}. A black-box description of the method goes as follows: given a positive integer $ k $, a Kalmar function $ f ( \vec{n} ) $, a simple-in-$ ( x_1 , \dots , x_k ) $ exponential polynomial $ P ( \vec{n} , x_1 , \dots , x_k ) $ and two arithmetic terms $ t ( \vec{n} ) > 1 $ and $ w ( \vec{n} ) $ such that $ f ( \vec{n} )  $ is equal to the cardinality of the set $ \{ \vec{a} \in \{ 0 , \dots , t ( \vec{n} ) - 1 \}^k : P ( \vec{n} , \vec{a} ) = 0 \} $ and $ P ( \vec{n} , \vec{a} ) $ belongs to the set $ \{ 0 , \dots , 2^{w ( \vec{n} )} - 1 \} $ for every point $ \vec{a} $ in $ \{ 0 , \dots , t ( \vec{n} ) - 1 \}^k $, the method computes an arithmetic term representing $ f ( \vec{n} ) $.

Given two integers $ a $ and $ w $ such that $ 0 \leq a < 2^w $, let $ \delta ( a , w ) $ denote the number $ ( 2^w - 1 ) ( 2^w - a + 1 ) $, which is equal to $ 2^{2 w} - 2^w a + a - 1 $.

The method is based on the following fundamental fact proven in both Marchenkov \cite[Lemma 6]{Marchenkov} and Mazzanti \cite[Lemma 4.5]{Mazzanti}.

\begin{lemma}\label{LemmaHW} Given two integers $ a $ and $ w $ such that $ 0 \leq a < 2^w $, we have that $$\HW ( \delta ( a , w ) ) = \begin{cases} 2w, & a = 0, \\ w, & a \neq 0.\end{cases}$$ \end{lemma}

Lemma \ref{LemmaSumToProd} will be of great utility in the proof of Lemma \ref{LemmaMazMar}.

\begin{lemma}\label{LemmaSumToProd} Given a positive integer $ k $, $ k $ non-negative integers $ u_1 , \dots , u_k $ and $ k + 1 $ integers $ v_1 , \dots , v_k , t $ exceeding one, we have that $$ \sum_{\vec{a} \in \{ 0 , \dots , t - 1 \}^k} \left( a_1^{u_1} v_1^{a_1} \dots a_k^{u_k} v_k^{a_k} \right) = G_{u_1} ( v_1 , t - 1 ) \dots G_{u_k} ( v_k , t - 1 ) . $$ \end{lemma}

{\bf Proof} Indeed, $$ \sum_{\vec{a} \in \{ 0 , \dots , t - 1 \}^k} \left( a_1^{u_1} v_1^{a_1} \dots a_k^{u_k} v_k^{a_k} \right) = \sum_{a_1 = 0}^{t-1} ( a_1^{u_1} v_1^{a_1} ) \dots \sum_{a_k = 0}^{t-1} ( a_k^{u_k} v_k^{a_k} ) = G_{u_1} ( v_1 , t - 1 ) \dots G_{u_k} ( v_k , t - 1 ) . $$ \qed

We express the main technique in the Lemma \ref{LemmaMazMar}, which will be applied in order to prove Theorem \ref{ThmMain}.

\begin{lemma}\label{LemmaMazMar} If $ P ( \vec{n} , \vec{x} ) $ is a simple-in-$ \vec{x} $ exponential polynomial, $ k $ is the (positive) length of the tuple $ \vec{x} $ and $ t ( \vec{n} ) > 1 $ and $ w ( \vec{n} ) $ are two arithmetic terms such that $ P ( \vec{n} , \vec{a} ) \in \{ 0 , \dots , 2^{w ( \vec{n} )} - 1 \} $ for every point $ \vec{a} \in \{ 0 , \dots , t ( \vec{n} ) - 1 \}^k $, then there is an arithmetic term that represents the cardinality of the set $ \{ \vec{a} \in \{ 0 , \dots , t ( \vec{n} ) - 1 \}^k : P ( \vec{n} , \vec{a} ) = 0 \} $. \end{lemma}

{\bf Proof} Given an arithmetic term $ t ( \vec{n} ) > 1 $ and a positive integer $ k $, let $ v $ denote the function that maps each point $ \vec{a} \in \{ 0 , \dots , t ( \vec{n} ) - 1 \}^k $ into the arithmetic term $ a_1 + a_2 t ( \vec{n} ) + \dots + a_k t ( \vec{n} )^{k - 1} $.

Observe that $ v $ enumerates the points of $ \{ 0 , \dots , t ( \vec{n} ) - 1 \}^k $ from $ 0 $ to $ t ( \vec{n} )^k - 1 $.

Let $$ M ( \vec{n} ) = \sum_{\vec{a} \in \{ 0 , \dots , t ( \vec{n} ) - 1 \}^k} \left( 2^{2 w ( \vec{n} ) v ( \vec{a} )} \delta ( P ( \vec{n} , \vec{a} ) , w ( \vec{n} ) ) \right) , $$ which is well-defined because $ 0 \leq P ( \vec{n} , \vec{a} ) < 2^{w ( \vec{n} )} $ for every $ \vec{a} \in \{ 0 , \dots , t ( \vec{n} ) - 1 \}^k $, and let $ d ( \vec{n} ) $ denote the cardinality of the set $ \{ \vec{a} \in \{ 0 , \dots , t ( \vec{n} ) - 1 \}^k : P ( \vec{n} , \vec{a} ) = 0 \} $.

Note that the binary representation of $ M ( \vec{n} ) $ is a concatenation of the binary representations of the $ t ( \vec{n} )^k  $ numbers $ \delta ( P ( \vec{n} , \vec{a} ) , w ( \vec{n} ) ) $ (with some extra zeros), which, by applying Lemma \ref{LemmaHW}, have at most $ 2 w ( \vec{n} ) $ ones each. Hence we have that $ \HW ( M ( \vec{n} ) ) = $ $$ \HW \left( \sum_{\vec{a} \in \{ 0 , \dots , t ( \vec{n} ) - 1 \}^k} \left( 2^{2 w ( \vec{n} ) v ( \vec{a} )} \delta ( P ( \vec{n} , \vec{a} ) , w ( \vec{n} ) ) \right) \right) = $$ $$ \sum_{\vec{a} \in \{ 0 , \dots , t ( \vec{n} ) - 1 \}^k} ( \HW ( \delta ( P ( \vec{n} , \vec{a} ) , w ( \vec{n} ) ) ) ) = $$ $$ d ( \vec{n} ) 2 w ( \vec{n} ) + ( t ( \vec{n} )^k - d ( \vec{n} ) ) w ( \vec{n} ) , $$ from which follows that $ \HW ( M ( \vec{n} ) ) / w ( \vec{n} ) - t ( \vec{n} )^k = d ( \vec{n} ) $.

Because $ P ( \vec{n} , \vec{x} ) $ is a simple-in-$ \vec{x} $ exponential polynomial, we know that, for some integer $ r \geq 1 $, it is a sum of $ r $ simple-in-$ \vec{x} $ exponential monomials $ m_1 ( \vec{n} , \vec{x} ) $, $ \dots $, $ m_r ( \vec{n} , \vec{x} ) $, each one containing at least one occurrence of a variable in $ \{ x_1 , \dots , x_k \} $, plus a (possibly zero) simple-in-$ \vec{x} $ exponential monomial $ \varepsilon ( \vec{n} ) $ (which is of the special kind of simple-in-$ \vec{x} $ exponential monomial in which, when written as in Expression \ref{ExprMonomial}, $ \beta_1 ( \vec{n} ) = \dots = \beta_k ( \vec{n} ) = 0 = \gamma_1 = \dots = \gamma_k $).

It only remains to express $ M ( \vec{n} ) $ as an arithmetic term: $$ M ( \vec{n} ) = \sum_{\vec{a} \in \{ 0 , \dots , t ( \vec{n} ) - 1 \}^k} \left( 2^{2 w ( \vec{n} ) v ( \vec{a} )} \delta ( P ( \vec{n} , \vec{a} ) , w ( \vec{n} ) ) \right) = $$ $$ \sum_{\vec{a} \in \{ 0 , \dots , t ( \vec{n} ) - 1 \}^k} \left( 2^{2 w ( \vec{n} ) v ( \vec{a} )} \delta \left( \varepsilon ( \vec{n} ) + m_1 ( \vec{n} , \vec{x} ) + \dots + m_r ( \vec{n} , \vec{x} ) , w ( \vec{n} ) \right) \right) = $$ $$ \sum_{\vec{a} \in \{ 0 , \dots , t ( \vec{n} ) - 1 \}^k} \left( 2^{2 w ( \vec{n} ) v ( \vec{a} )} \left( 2^{w ( \vec{n} )} - 1 \right) \left( 2^{w ( \vec{n} )} - \varepsilon ( \vec{n} ) - m_1 ( \vec{n} , \vec{x} ) - \dots - m_r ( \vec{n} , \vec{x} ) + 1 \right) \right) = $$ $$ \sum_{\vec{a} \in \{ 0 , \dots , t ( \vec{n} ) - 1 \}^k} \left( 2^{2 w ( \vec{n} ) v ( \vec{a} )} \left( 2^{w ( \vec{n} )} - 1 \right) \left( 2^{w ( \vec{n} )} - \varepsilon ( \vec{n} ) + 1 \right) \right) + $$ $$ \sum_{j = 1}^r \left( \sum_{\vec{a} \in \{ 0 , \dots , t ( \vec{n} ) - 1 \}^k} \left( 2^{2 w ( \vec{n} ) v ( \vec{a} )} \left( 2^{w ( \vec{n} )} - 1 \right) ( - m_j ( \vec{n} , \vec{x} ) ) \right) \right) . $$

Now we study each simple-in-$ \vec{x} $ exponential monomial separately.

First, by applying Lemma \ref{LemmaSumToProd}, we have that $$ \sum_{\vec{a} \in \{ 0 , \dots , t ( \vec{n} ) - 1 \}^k} \left( 2^{2 w ( \vec{n} ) v ( \vec{a} )} \left( 2^{w ( \vec{n} )} - 1 \right) \left( 2^{w ( \vec{n} )} - \varepsilon ( \vec{n} ) + 1 \right) \right) = $$ $$ \left( 2^{w ( \vec{n} )} - 1 \right) \left( 2^{w ( \vec{n} )} - \varepsilon ( \vec{n} ) + 1 \right) \sum_{\vec{a} \in \{ 0 , \dots , t ( \vec{n} ) - 1 \}^k} \left( 2^{2 w ( \vec{n} ) v ( \vec{a} )} \right) = $$ $$ \left( 2^{w ( \vec{n} )} - 1 \right) \left( 2^{w ( \vec{n} )} - \varepsilon ( \vec{n} ) + 1 \right) \sum_{\vec{a} \in \{ 0 , \dots , t ( \vec{n} ) - 1 \}^k} \left( \left( 2^{2 w ( \vec{n} )} \right)^{a_1} \dots \left( 2^{2 w ( \vec{n} ) {t ( \vec{n} )}^{k - 1}} \right)^{a_k} \right) = $$ $$ \left( 2^{w ( \vec{n} )} - 1 \right) \left( 2^{w ( \vec{n} )} - \varepsilon ( \vec{n} ) + 1 \right) G_0 \left( 2^{2 w ( \vec{n} )}, t ( \vec{n} ) - 1 \right) \dots G_0 \left( 2^{2 w ( \vec{n} ) {t ( \vec{n} )}^{k - 1}}, t ( \vec{n} ) - 1 \right) = $$ $$ \left( 2^{w ( \vec{n} )} - 1 \right) \left( 2^{w ( \vec{n} )} - \varepsilon ( \vec{n} ) + 1 \right) \left( 2^{2 w ( \vec{n} ) {t ( \vec{n} )}^k } - 1 \right) / \left( 2^{2 w ( \vec{n} )} - 1 \right) = $$ \begin{equation}\label{TermC} \left( 2^{w ( \vec{n} )} - \varepsilon ( \vec{n} ) + 1 \right) \left( 2^{2 w ( \vec{n} ) {t ( \vec{n} )}^k } - 1 \right) / \left( 2^{w ( \vec{n} )} + 1 \right) , \end{equation} as many numerators and denominators cancel with each other (recall Identity \ref{EqGeomSum}).

Now, let $ m ( \vec{n} , \vec{x} ) \in \{ m_1 ( \vec{n} , \vec{x} ) , \dots , m_r ( \vec{n} , \vec{x} ) \} $.

Then, again by applying Lemma \ref{LemmaSumToProd}, and writing $ m ( \vec{n} , \vec{x} ) $ as in Expression \ref{ExprMonomial}: $$ \sum_{\vec{a} \in \{ 0 , \dots , t ( \vec{n} ) - 1 \}^k} \left( 2^{2 w ( \vec{n} ) v ( \vec{a} )} \left( 2^{w ( \vec{n} )} - 1 \right) ( - m ( \vec{n} , \vec{x} ) ) \right) = $$ $$ - \left( 2^{w ( \vec{n} )} - 1 \right) \sum_{\vec{a} \in \{ 0 , \dots , t ( \vec{n} ) - 1 \}^k} \left( 2^{2 w ( \vec{n} ) v ( \vec{a} )} m ( \vec{n} , \vec{a} ) \right) = $$ $$ - \left( 2^{w ( \vec{n} )} - 1 \right) \sum_{\vec{a} \in \{ 0 , \dots , t ( \vec{n} ) - 1 \}^k} \left( 2^{2 w ( \vec{n} ) v ( \vec{a} )} \alpha ( \vec{n} ) b_1 ( \vec{n} )^{\beta_1 ( \vec{n} ) a_1} \dots b_k ( \vec{n} )^{\beta_k ( \vec{n} ) a_k} a_1^{\gamma_1} \dots a_k^{\gamma_k} \right) = $$ $$ - \left( 2^{w ( \vec{n} )} - 1 \right) \alpha ( \vec{n} ) \sum_{\vec{a} \in \{ 0 , \dots , t ( \vec{n} ) - 1 \}^k} \left( \left( 2^{2 w ( \vec{n} )} b_1 ( \vec{n} )^{\beta_1 ( \vec{n} )} \right)^{a_1} a_1^{\gamma_1} \dots \left( 2^{2 w ( \vec{n} ) {t ( \vec{n} )}^{k - 1}} b_k ( \vec{n} )^{\beta_k ( \vec{n} )} \right)^{a_k} a_k^{\gamma_k} \right) = $$ \begin{equation}\label{TermA} - \left( 2^{w ( \vec{n} )} - 1 \right) \alpha ( \vec{n} ) G_{\gamma_1} \left( 2^{ 2 w ( \vec{n} )} b_1 ( \vec{n} )^{\beta_1 ( \vec{n} )}, t ( \vec{n}) - 1 \right) \dots G_{\gamma_k} \left( 2^{ 2 w ( \vec{n} ) {t ( \vec{n} )}^{k - 1}} b_k ( \vec{n} )^{\beta_k ( \vec{n} )}, t ( \vec{n} ) -1 \right) . \end{equation} \qed

The expressions \ref{TermC} and \ref{TermA} from the proof of Lemma \ref{LemmaMazMar} will be denoted, respectively, by $ \mathcal{C} ( \varepsilon ( \vec{n} ) , k ) $ and $ \mathcal{A} ( m ( \vec{n} , \vec{x} ) , k ) $. It is important to remark that, even in the case that $ \varepsilon ( \vec{n} ) $ is zero, the expression $ \mathcal{C} ( \varepsilon ( \vec{n} ) , k ) $ is non-zero.

In the following sections, we will often consider the bound $ w ( \vec{n} ) $ from the statement of Lemma \ref{LemmaMazMar} larger than necessary, in order to keep the proofs relatively simple. But the reader should keep in mind that it could be sharpened.

\begin{theo}\label{ThmMain} If $ P ( \vec{n} , \vec{x} ) $ is an algebraic sum of arithmetic terms, $ k $ is the (positive) length of $ \vec{x} $ and $ t ( \vec{n} ) $ is an arithmetic term exceeding one, then there is an arithmetic term in variables $ \vec{n} $, and built up by using the Hamming weight and generalized geometric progressions, that represents the cardinality of the set $ \{ \vec{a} \in \{ 0 , \dots , t ( \vec{n} ) - 1 \}^k : P ( \vec{n} , \vec{a} ) = 0 \} $. \end{theo}

{\bf Proof} First, add new variables $ \vec{y} $ with which to encode all the necessary subterms until having a sum of squares of simple-in-$ ( \vec{x} , \vec{y} ) $ exponential polynomials whose expansion $ Q ( \vec{n} , \vec{x} , \vec{y} ) $ is a simple-in-$ ( \vec{x} , \vec{y} ) $ exponential polynomial such that $$ \forall \ \vec{r} , \vec{a} , \vec{b} \in \N \ [ Q ( \vec{r} , \vec{a} , \vec{b} ) = 0 \ \Rightarrow \ P ( \vec{r} , \vec{a} ) = 0 ] $$ and \begin{equation}\label{UniquenessCondition} \forall \ \vec{r} , \vec{a} \in \N \ [ P ( \vec{r} , \vec{a} ) = 0 \ \Rightarrow \ \exists ! \ \vec{b} \in \N \ Q ( \vec{r} , \vec{a} , \vec{b} ) = 0 ] \end{equation} (see an example below).

Let $ f $ be the length of the tuple $ \vec{y} $.

Then there is some arithmetic term $ \theta ( \vec{n} ) $ such that, if $ \vec{a} \in \{ 0 , \dots , t ( \vec{n} ) - 1 \}^k $ and $  P ( \vec{n} , \vec{a} ) = 0 $, then $ \vec{b} \in \{ 0 , \dots , \theta ( \vec{n} ) - 1 \}^f $ for the corresponding solution $ ( \vec{a} , \vec{b} ) $ of the equation $ Q ( \vec{n} , \vec{x} , \vec{y} ) = 0 $.

Now, find an arithmetic term $ w ( \vec{n} ) $ such that the inequality $ Q ( \vec{n} , \vec{a} , \vec{b} ) < 2^{w ( \vec{n} )} $ holds for every $ ( \vec{a} , \vec{b} ) \in \{ 0 , \dots , \max ( t ( \vec{n} ) , \theta ( \vec{n} ) ) - 1 \}^{k + f} $. And we know that $ Q ( \vec{n} , \vec{x} , \vec{y} ) \geq 0 $ because it equals a sum of squares.

By using the construction from the proof of Lemma \ref{LemmaMazMar} for $ Q ( \vec{n} , \vec{x} , \vec{y} ) $ and the arithmetic terms $ \max ( t ( \vec{n} ) , \theta ( \vec{n} ) ) $ and $ w ( \vec{n} ) $, obtain an arithmetic term $ d ( \vec{n} ) $ representing the cardinality of the set $ \{ ( \vec{a} , \vec{b} ) \in \{ 0 , \dots , \max ( t ( \vec{n} ) , \theta ( \vec{n} ) ) - 1 \}^{k + f} : Q ( \vec{n} , \vec{a} , \vec{b} ) = 0 \} $.

Because of Condition \ref{UniquenessCondition}, we can conclude that $ d ( \vec{n} ) $ also represents the cardinality of the set $ \{ \vec{a} \in \{ 0 , \dots , t ( \vec{n} ) - 1 \}^k : P ( \vec{n} , \vec{a} ) = 0 \} $. \qed

For example, if in $ P ( \vec{n} , \vec{x} ) $ from Theorem \ref{ThmMain} the arithmetic term $ x_1^{x_2} $ occurs, then we have, by applying Identity \ref{EqMarchenkov}, that: $$ P ( \vec{n} , \vec{x} ) = 0 \ \Leftrightarrow $$ $$ [ y_1 = x_1^{x_2} \ \wedge \ P ( \vec{n} , \vec{x} , y_1 ) = 0 ] \ \Leftrightarrow $$ $$ [ y_1 = 2^{( x_1 x_2 + x_1 + 1 ) x_2} \bmod ( 2^{x_1 x_2 + x_1 + 1} - x_1 ) \ \wedge \ P ( \vec{n} , \vec{x} , y_1 ) = 0 ] \ \Leftrightarrow $$ $$ [ y_2 = x_1 x_2 + x_1 + 1 \ \wedge \ y_3 = y_2 x_2 \ \wedge \ y_1 = 2^{y_3} \bmod ( 2^{y_2} - x_1 ) \ \wedge \ P ( \vec{n} , \vec{x} , y_1 ) = 0 ] \ \Leftrightarrow $$ $$ [ y_2 = x_1 x_2 + x_1 + 1 \ \wedge \ y_3 = y_2 x_2 \ \wedge \ ( 2^{y_2} - x_1 ) y_4 + y_1 = 2^{y_3} \wedge y_1 < 2^{y_2} - x_1 \ \wedge \ P ( \vec{n} , \vec{x} , y_1 ) = 0 ] \ \Leftrightarrow $$ $$ [ y_2 = x_1 x_2 + x_1 + 1 \ \wedge \ y_3 = y_2 x_2 \ \wedge \ ( 2^{y_2} - x_1 ) y_4 + y_1 = 2^{y_3} \ \wedge \ y_1 + y_5 + 1 = 2^{y_2} - x_1 \ \wedge \ P ( \vec{n} , \vec{x} , y_1 ) = 0 ] \ \Leftrightarrow $$ $$ ( y_2 - x_1 x_2 - x_1 - 1 )^2 + ( y_3 - y_2 x_2 )^2 + ( ( 2^{y_2} - x_1 ) y_4 + y_1 - 2^{y_3} )^2 + ( y_1 + y_5 + 1 - 2^{y_2} - x_1 )^2 + P ( \vec{n} , \vec{x} , y_1 )^2 = 0 \ \Leftrightarrow $$ $$ Q ( \vec{n} , \vec{x} , y_1 , \dots , y_5 ) = 0 . $$

We say that the \textbf{non-exponential occurrences} of an expression $ t $ in (the expansion of) a simple-in-$ \vec{x} $ exponential polynomial $ P ( \vec{n} , \vec{x} ) $ are those in which $ t $ appears as a factor of some simple-in-$ \vec{x} $ exponential mononomial of $ P ( \vec{n} , \vec{x} ) $.

For example, the non-exponential occurrences of $ x_1^2 $ in $ x_1^2 2^{x_1^2} - x_1^2 $ are the first and the third one, and the non-exponential occurrence of $ x_1 $ in Expression \ref{ExprMonomial} is the second one.

\begin{lemma}\label{LemmaExp012} If $ P ( \vec{n} , \vec{x} ) $ is a simple-in-$ \vec{x} $ exponential polynomial, then there is a simple-in-$ ( \vec{x} , \vec{y} ) $ exponential polynomial $ Q ( \vec{n} , \vec{x} , \vec{y} ) \geq 0 $ that satisfies the following conditions.
\begin{enumerate}
    \item No non-exponential occurrence in $ Q ( \vec{n} , \vec{x} , \vec{y} ) $ of a variable in $ ( \vec{x} , \vec{y} ) $ has an exponent larger than two.
    \item $ \forall \ \vec{r} , \vec{a} , \vec{b} \in \N \ [ Q ( \vec{r} , \vec{a} , \vec{b} ) = 0 \ \Rightarrow \ P ( \vec{r} , \vec{a} ) = 0 ] $.
    \item $ \forall \ \vec{r} , \vec{a} \in \N \ [ P ( \vec{r} , \vec{a} ) = 0 \ \Rightarrow \ \exists ! \ \vec{b} \in \N \ Q ( \vec{r} , \vec{a} , \vec{b} ) = 0 ] $.
\end{enumerate} \end{lemma}

{\bf Proof} Let $ k $ be the (positive) length of the tuple $ \vec{x} $, and let $ h $ be the largest (positive) exponent from among all the non-exponential occurrences in $ P ( \vec{n} , \vec{x} ) $ of a variable $ x $ in $ \vec{x} $.

We consider new variables $ y_1, \dots, y_h $ and the following $ h $ equations: \begin{eqnarray*} x - y_1 & = & 0 , \\ y_1 x - y_2 & = & 0 , \\ \vdots & \vdots & \vdots \\ y_{h - 1} x - y_h & = & 0 . \end{eqnarray*}

For every $ i \in \{ 1 , \dots , h \} $, we replace in $ P ( \vec{n} , \vec{x} ) $ all the non-exponential occurrences of $ x^i $ with the corresponding variable $ y_i $.

And the same procedure is done for all variables $ x_1 , \dots , x_k $ (of course, by adding new variables of the form $ y_i $ each time).

In the end, we get a simple-in-$ ( \vec{x} , \vec{y} ) $ exponential polynomial $ P ( \vec{n} , \vec{x} , \vec{y} ) $ in which the variables in $ \vec{x} $ have no non-exponential occurrences, and in which variables in $ \vec{y} $ have only non-exponential occurrences with exponent one.

There are two kinds of simple-in-$ ( \vec{x} , \vec{y} ) $ exponential monomials in the expansion of $ P (\vec{n} , \vec{x} , \vec{y} )^2 $: those of the form $ m^2 $ and those of the form $ 2 m m ' $, where $ m $ and $ m ' $ are simple-in-$ ( \vec{x} , \vec{y} ) $ exponential monomials of $ P (\vec{n} , \vec{x} , \vec{y} ) $. In $ m^2 $, the non-exponential occurrences of the variables in $ \vec{y} $ have exponent two. In $ 2 m m ' $, the non-exponential occurrences of the variables in $ \vec{y} $ which are common to both $ m $ and $ m ' $ have exponent two; and the non-exponential occurrences of all the other variables in $ \vec{y} $ have exponent one. Therefore, no non-exponential occurrence in the expansion of $ P ( \vec{n} , \vec{x} , \vec{y} )^2 $ of a variable in $ \vec{y} $ has an exponent larger than two.

Finally, if $ S_1 ( \vec{n} , \vec{x} , \vec{y} ) $, $ \dots $, $ S_f ( \vec{n} , \vec{x} , \vec{y} ) $ are the left-hand sides of the introduced equations, then we define $ Q ( \vec{n} , \vec{x} , \vec{y} ) $ as the expansion of $ S_1 ( \vec{n} , \vec{x} , \vec{y} )^2 + \dots + S_f ( \vec{n} , \vec{x} , \vec{y} )^2 + P ( \vec{n} , \vec{x} , \vec{y} )^2 $ and all the conditions from the statement become clear. \qed

As we advanced in Section \ref{SectionUsefulTerms}, Corollary \ref{CorFourArithmeticTerms} shows that, together with the Hamming weight, only the first three kinds of generalized geometric progressions are necessary to produce arithmetic terms. However, some functions might have much shorter formulas when expressed in terms of generalized geometric progressions of higher kind.

\begin{cor}\label{CorFourArithmeticTerms} If $ P ( \vec{n} , \vec{x} ) $ is an algebraic sum of arithmetic terms, $ k $ is the (positive) length of the tuple $ \vec{x} $ and $ t ( \vec{n} ) $ is an arithmetic term exceeding one, then there is an arithmetic term in variables $ \vec{n} $, and built up by using the Hamming weight and generalized geometric progressions of up to the second kind, that represents the cardinality of the set $ \{ \vec{a} \in \{ 0 , \dots , t ( \vec{n} ) - 1 \}^k : P ( \vec{n} , \vec{a} ) = 0 \} $. \end{cor}

{\bf Proof} As in the proof of Theorem \ref{ThmMain}, we add new variables $ \vec{y} $ with which to encode all the necessary subterms until having a sum of squares of simple-in-$ ( \vec{x} , \vec{y} ) $ exponential polynomials whose expansion $ Q ( \vec{n} , \vec{x} , \vec{y} ) $ is a simple-in-$ ( \vec{x} , \vec{y} ) $ exponential polynomial such that $ \forall \ \vec{r} , \vec{a} , \vec{b} \in \N \ [ Q ( \vec{r} , \vec{a} , \vec{b} ) = 0 \ \Rightarrow \ P ( \vec{r} , \vec{a} ) = 0 ] $ and $ \forall \ \vec{r} , \vec{a} \in \N \ [ P ( \vec{r} , \vec{a} ) = 0 \ \Rightarrow \ \exists ! \ \vec{b} \in \N \ Q ( \vec{r} , \vec{a} , \vec{b} ) = 0 ] $.

We then apply Lemma \ref{LemmaExp012} in order to get a  simple-in-$ ( \vec{x} , \vec{y} , \vec{z} ) $ exponential polynomial $ R ( \vec{n} , \vec{x} , \vec{y} , \vec{z} ) \geq 0 $ that satisfies the following conditions.
\begin{enumerate}
    \item No non-exponential occurrence in $ R ( \vec{n} , \vec{x} , \vec{y} , \vec{z} ) $ of a variable in $ ( \vec{x} , \vec{y} , \vec{z} ) $ has an exponent larger than two.
    \item $ \forall \ \vec{r} , \vec{a} , \vec{b} , \vec{c} \in \N \ [ R ( \vec{r} , \vec{a} , \vec{b} , \vec{c} ) = 0 \ \Rightarrow \ Q ( \vec{r} , \vec{a} , \vec{b} ) = 0 ] $.
    \item $ \forall \ \vec{r} , \vec{a} , \vec{b} \in \N \ [ Q ( \vec{r} , \vec{a} , \vec{b} ) = 0 \ \Rightarrow \ \exists ! \ \vec{c} \in \N \ R ( \vec{r} , \vec{a} , \vec{b} , \vec{c} ) = 0 ] $.
\end{enumerate}

Now, let $ f $ be the length of the tuple $ ( \vec{y} , \vec{z} ) $.

Then there is some arithmetic term $ \theta ( \vec{n} ) $ such that, if $ \vec{a} \in \{ 0 , \dots , t ( \vec{n} ) - 1 \}^k $ and $ P ( \vec{n} , \vec{a} ) = 0 $, then $ ( \vec{b} , \vec{c} ) \in \{ 0 , \dots , \theta ( \vec{n} ) - 1 \}^f $ for the corresponding solution $ ( \vec{a} , \vec{b} , \vec{c} ) $ of the equation $ R ( \vec{n} , \vec{x} , \vec{y} , \vec{z} ) = 0 $.

We find an arithmetic term $ w ( \vec{n} ) $ such that the inequality $ R ( \vec{n} , \vec{a} , \vec{b} , \vec{c} ) < 2^{w ( \vec{n} )} $ holds for every $( \vec{a} , \vec{b} , \vec{c} ) \in \{ 0 , \dots , \max ( t ( \vec{n} ) , \theta ( \vec{n} ) ) - 1 \}^{k + f} $.

And finally, by using the construction from the proof of Lemma \ref{LemmaMazMar} for $ R ( \vec{n} , \vec{x} , \vec{y} , \vec{z} ) $ and the arithmetic terms $ \max ( t ( \vec{n} ) , \theta ( \vec{n} ) ) $ and $ w ( \vec{n} ) $, we find an arithmetic term representing the cardinality of the set $ \{ \vec{a} \in \{ 0 , \dots , t ( \vec{n} ) - 1 \}^k : P ( \vec{n} , \vec{a} ) = 0 \} $. As all the non-exponential occurrences in $ R ( \vec{n} , \vec{x} , \vec{y} , \vec{z} ) $ of the variables in $ ( \vec{x} , \vec{y} , \vec{z} ) $ have exponent one or two, only the generalized geometric progressions $ G_0 $, $ G_1 $ and $ G_2 $ are necessary to build this arithmetic term (recall the proof of Lemma \ref{LemmaMazMar}). \qed

\section{The number-of-divisors function}\label{SectionTau}

\begin{lemma}\label{LemmaSetTau} If $ n $ is a positive integer, then $ \tau ( n ) $ is equal to the cardinality of the set $$ \{ ( a , b ) \in \{ 0 , \dots , n \}^2 : n - a b = 0 \} . $$ \end{lemma}

{\bf Proof} We have that $ \tau ( n ) = $

$ | \{ a \in \{ 1 , \dots , n \} : \textrm{$ a $ divides $ n $} \} | = $

$ | \{ a \in \{ 1 , \dots , n \} : \textrm{exists $ b \in \{ 1 , \dots , n \} $ such that $ a b = n $} \} | = $

$ | \{ ( a , b ) \in \{ 0 , \dots , n \}^2 : n - ab = 0 \} | $. \qed

\begin{lemma}\label{LemmaMajTau} If $ n $ is a positive integer and $ ( a , b ) \in \{ 0 , \dots , n \}^2$, then $ ( n - a b )^2 < 2^{n + 4} $. \end{lemma}

{\bf Proof} The largest number of the form $ | n - a b | $, where $ ( a , b ) \in \{ 0 , \dots , n \}^2 $, is clearly $ n^2 - n $.

And it is easy to see that $ ( n^2 - n )^2 < 2^{n + 4} $. \qed

\begin{theo}\label{ThmTau} The function $ \tau ( n ) $ (for positive integer arguments $ n $) can be represented by the arithmetic term $$ \HW(M(n))/(n+4) - (n+1)^2,$$ where $ M ( n ) $ is equal to \begin{equation}\label{ExprTau} \mathcal{C} ( n^2 , 2 ) + \mathcal{A} ( - 2 n x_1 x_2 , 2 ) + \mathcal{A} ( x_1^2 x_2^2 , 2 ) . \end{equation} \end{theo}

{\bf Proof} Let $ P ( n , x_1 , x_2 ) = ( n - x_1 x_2 )^2 $, $ t ( n ) = n + 1 $ and $ w ( n ) = n + 4 $.

By applying Lemma \ref{LemmaMajTau}, $ P ( n , a , b ) \in \{ 0 , \dots , 2^{w ( n )} - 1 \} $ for every point $ ( a , b ) \in \{ 0 , \dots , t ( n ) - 1 \}^2 $. Therefore, we can instantiate the proof of Lemma \ref{LemmaMazMar} to this particular case.

Notice that $$ P ( n , x_1 , x_2 ) = n^2 - 2 n x_1 x_2 + x_1^2 x_2^2 , $$ so we define $ M ( n ) $ as Expression \ref{ExprTau}.

Then we have that \begin{equation}\label{EqTau} | \{ ( a , b ) \in \{ 0 , \dots , t ( n ) - 1 \}^2 : P ( n , a , b ) = 0 \} | = \HW ( M ( n ) ) / w ( n ) - t ( n )^2 . \end{equation}

The left-hand side of Identity \ref{EqTau} is equal to $$ | \{ ( a , b ) \in \{ 0 , \dots , t ( n ) - 1 \}^2 : n - a b = 0 \} | $$ and hence, by applying Lemma \ref{LemmaSetTau}, to $ \tau ( n ) $. \qed

Exceptionally, in order to give an idea of the form of the expression $ M ( n ) $ from Theorem \ref{ThmTau}, we also represent it explicitly:

$ ( 2^{n + 4} - n^2 + 1 ) ( 2^{n + 4} + 1 )^{- 1} ( 2^{2 ( n + 4 ) ( n + 1 )^2} - 1 ) $

$ + 2^{2 ( n + 4 ) ( n + 2 ) + 1} n ( 2^{n + 4} - 1 ) ( 2^{2 ( n + 4 )} - 1 )^{- 2} ( 2^{2 ( n + 4 ) ( n + 1 )} - 1 )^{- 2} $

$ \cdot ( 2^{2 ( n + 4 ) ( n + 1 )} n - 2^{2 ( n + 4 ) n} ( n + 1 ) + 1 ) ( 2^{2 ( n + 4 ) ( n + 1 )^2} n - 2^{2 ( n + 4 ) ( n + 1 ) n} ( n + 1 ) + 1 ) $

$ - 2^{2 ( n + 4 ) ( n + 2 )} ( 2^{n + 4} - 1 ) ( 2^{2 ( n + 4 )} - 1 )^{- 3} ( 2^{2 ( n + 4 ) ( n + 1 )} - 1 )^{- 3} $

$ \cdot ( 2^{2 ( n + 4 ) ( n + 2 )} n^2 - 2^{2 ( n + 4 ) ( n + 1 )} ( 2 n^2 + 2 n - 1 ) + 2^{2 ( n + 4 ) n} ( n + 1 )^2 - 2^{2 ( n + 4 )} - 1 ) $

$ \cdot ( 2^{2 ( n + 4 ) ( n + 1 ) ( n + 2 )} n^2 - 2^{2 ( n + 4 ) ( n + 1 )^2} ( 2 n^2 + 2 n - 1 ) + 2^{2 ( n + 4 ) ( n + 1 ) n} ( n + 1 )^2 - 2^{2 ( n + 4 ) ( n + 1 )} - 1 ) $.

Observe that the proof of Theorem \ref{ThmTau} is a particularization of the proof of Lemma \ref{LemmaMazMar} to the concrete case of $ \tau $. In the next proofs we will omit such details, since the procedure is completely analogous.

\begin{cor} The set of prime numbers has the quantifier-free arithmetic-term representation $ \tau(n) = 2 $.\end{cor}

Proposition \ref{PropTauPadicVal} allows one to express the $ p $-adic valuation in terms of the number-of-divisors function. The case in which $ p = 2 $ is due to Stearns \& Yanev (see \href{https://oeis.org/A007814}{\texttt{OEIS A007814}}).

\begin{prop}\label{PropTauPadicVal} If $ p $ is a prime and $ n $ is a positive integer, then $$ \nu_p ( n ) ( \tau ( p n ) - \tau ( n ) ) = 2 \tau ( n ) - \tau ( p n ) . $$ \end{prop}

{\bf Proof} Let $ x = \nu_p ( n ) $, let $ y $ be the positive integer such that $ p^x y = n $ and let $ z = \tau ( y ) $.

Note that $ \gcd ( p , y ) = 1 $.

Then $ \tau ( n ) = \tau ( p^x y ) = ( x + 1 ) z $ and $ \tau ( p n ) = \tau ( p^{x + 1} y ) = ( x + 2 ) z $.

The identity to prove becomes $ x ( ( x + 2 ) z - ( x + 1 ) z ) = 2 ( x + 1) z - ( x + 2 ) z $, which is clearly true. \qed

\section{The sum-of-divisors function}\label{SectionSigma}

\begin{lemma}\label{LemmaSetSigma} If $ n $ is a positive integer, then $ \sigma ( n ) $ is equal to the cardinality of the set $$ \{ ( a , b , c ) \in \{ 0 , \dots , n \}^3 : n - ( a + b + 1 ) c = 0 \} . $$ \end{lemma}

{\bf Proof} Let $ A $ be the set of the statement and let $ D ( n ) $ be the set of divisors of $ n $.

For every $ x \in D ( n ) $, let $ B_x $ be the set $ \{ ( a , x - a - 1 , n / x ) : a \in \{ 0 , \dots , n \}^3 \ \textrm{and} \ x - a - 1 \geq 0 \} . $

For every $ x \in D ( n ) $ and every $ a \in \{ 0 , \dots , n \} $, the condition $ 0 \leq x - a - 1 $ implies that $ a \leq x - 1 $. Thus, the cardinality of $ B_x $ is $ x $.

Let $ B $ be the set $ \bigcup ( \{ B_x : x \in D ( n ) \} ) $.

The elements of $ \{ B_x : x \in D ( n ) \} $ are obviously pairwise disjoint, from which follows that $$ | B | = \sum_{x \in D ( n )} ( | B_x | ) = \sum_{x \in D ( n )} ( x ) = \sigma ( n ) . $$

The inclusion $ B \subseteq A $ is immediate, so it only remains to check that $ A \subseteq B $. Let $ ( a , b , c ) \in A $ and let $ d = a + b + 1 $.

Then it is clear that $ d | n $ and $ ( a , b , c ) = ( a , d - a - 1 , n / d ) \in B_d \subseteq B $. \qed

\begin{lemma}\label{LemmaMajSigma} If $ n $ is a positive integer and $ ( a , b , c ) \in \{ 0 , \dots , n \}^3 $, then $ ( n - ( a + b + 1 ) c )^2 < 2^{n + 7} $. \end{lemma}

{\bf Proof} The largest number of the form $ | n - ( a + b + 1 ) c | $, where $ ( a , b , c ) \in \{ 0 , \dots , n \}^3 $, is clearly $ 2 n^2 $. And it is easy to see that $ ( 2 n^2 )^2 < 2^{n + 7} $. \qed

\begin{theo}\label{ThmSigma} The function $ \sigma(n) $ (for positive integer arguments $ n $) can be represented by the arithmetic term $$ \HW(M(n))/(n+7) - (n+1)^3 , $$ where $ M ( n ) $ is equal to $$ \begin{array}{llll} \mathcal{C} ( n^2 , 3 ) & + \mathcal{A} ( x_1^2 x_3^2 , 3 ) & + \mathcal{A} ( - 2 n x_1 x_3 , 3 ) & + \mathcal{A} ( 2 x_1 x_2 x_3^2 , 3 ) \\ & + \mathcal{A} ( x_2^2 x_3^2 , 3 ) & + \mathcal{A} ( - 2 n x_2 x_3 , 3 ) & + \mathcal{A} ( 2 x_1 x_3^2 , 3 ) \\ & + \mathcal{A} ( x_3^2 , 3 ) & + \mathcal{A} ( - 2 n x_3 , 3 ) & + \mathcal{A} ( 2 x_2 x_3^2 , 3 ) . \end{array} $$ \end{theo}

{\bf Proof} It follows from Lemma \ref{LemmaSetSigma} and Lemma \ref{LemmaMajSigma} by emulating the proof of Lemma \ref{LemmaMazMar}. \qed

A number $ n $ is said to be \textbf{perfect} if, and only if, $ \sigma(n) = 2 n $ (cf.\ K\v{r}\'{i}\v{z}ek et al.\ \cite[p.\ 85]{KrizekEtAl}).

\begin{cor} The set of perfect numbers has a quantifier-free arithmetic-term representation. \end{cor}

\section{Euler's totient function}\label{SectionPhi}

\begin{lemma}\label{LemmaSetPhi} If $ n $ is an integer exceeding one, then $ \varphi ( n ) $ is equal to the cardinality of the set $$ \{ ( a , b , c ) \in \{ 0 , \dots , n \}^3 : ab - cn - 1 = 0 \} . $$ \end{lemma}

{\bf Proof} We have that $ \varphi ( n ) = $

$ | \{ a \in \{ 1 , \dots , n - 1 \} : \gcd ( a , n ) = 1 \} | = $

$ | \{ a \in \{ 1 , \dots , n - 1 \} : \textrm{exists $ b \in \{  1 , \dots , n - 1 \} $ such that $ b = \inv ( a , n ) $} \} | = $

$ | \{ a \in \{ 1 , \dots , n - 1 \} : \textrm{exists $ b \in \{  1 , \dots , n - 1 \} $ such that $ a b \equiv 1 $ (mod $ n $)} \} | = $

$ | \{ a \in \{ 1 , \dots , n - 1 \} : \textrm{exist $ b , c \in \{ 1 , \dots , n - 1 \} $ such that $ a b - 1 = c n $} \} | = $

$ | \{ ( a , b , c ) \in \{ 0 , \dots , n \}^3 : ab - cn - 1 = 0 \} | $. \qed

\begin{lemma}\label{LemmaMajPhi} If $ n $ is a positive integer and $ ( a , b , c ) \in \{ 0 , \dots , n \}^3 $, then $ ( a b - c n - 1 )^2 < 2^{n + 5} $. \end{lemma}

{\bf Proof} The largest number of the form $ | ab - nc - 1| $, where $ ( a , b , c ) \in \{ 0 , \dots , n \}^3 $, is clearly $ n^2 + 1 $.
And it is easy to see that $ ( n^2 + 1 )^2 < 2^{n + 5} $. \qed

\begin{theo}\label{ThmPhi} The function $ \varphi ( n ) $ (for integer arguments $ n \geq 2 $) can be represented by the arithmetic term $$ \HW ( M ( n ) ) / ( n + 5 ) - ( n + 1 )^3 , $$ where $ M ( n ) $ is equal to $$ \begin{array}{lll} \mathcal{C} ( 1 , 3 ) & + \mathcal{A} ( x_1^2 x_2^2 , 3 ) & +  \mathcal{A} ( n^2 x_3^2 , 3 ) \\ + \mathcal{A} ( -2nx_1x_2x_3 , 3 ) & + \mathcal{A} ( -2x_1 x_2 , 3 ) & + \mathcal{A} ( 2 n x_3 , 3 ) . \end{array} $$ \end{theo}

{\bf Proof} It follows from Lemma \ref{LemmaSetPhi} and Lemma \ref{LemmaMajPhi} by emulating the proof of Lemma \ref{LemmaMazMar}. \qed

\section{The modular inverse}\label{SectionInv}

\begin{lemma}\label{LemmaSetInv} Given two coprime integers $ n \geq 2 $ and $ m \in \{ 1 , \dots , n - 1 \} $, the number $ \inv ( m , n ) $ is equal to the cardinality of the set $ \{ ( a , b , c , d ) \in \{ 0 , \dots , n \}^4 : ( m a - n b - 1 )^2 + ( a - c - d - 1 )^2 = 0 \} $. \end{lemma}

{\bf Proof} Let $ x $ be the number $ \inv ( m , n ) $, let $ A $ be the set of the statement and let $ B $ be the set $$ \{ ( x , ( m x - 1 ) / n , c , x - c - 1 ) : c \in \{0, \dots, n\} \ \textrm{and} \ x - c - 1 \geq 0 \} . $$

The number $ x - c - 1 $ is non-negative, so $ c \in \{ 0 , \dots , x - 1 \} $ and hence the cardinality of $ B $ is $ x $.

The inclusion $ B \subseteq A $ is immediate, so it only remains to check that $ A \subseteq B $. Let $ ( a , b , c , d ) \in A $. 

The conditions $ m a = n b + 1 $, $ \gcd ( m , n ) = 1 $, $ n \geq 2 $ and $ m \in \{ 1 , \dots , n - 1 \} $ imply that $ a = x $ and $ b = ( m x - 1 ) / n $. And the condition $ a - c - d - 1 = 0 $ then implies that $ d = x - c - 1 $ and $ c \leq x - 1 $, so we can conclude that $ (a, b, c, d) \in B $. \qed

\begin{lemma}\label{LemmaMajInv} If $ n $ is an integer exceeding one, $ m \in \{ 1 , \dots , n - 1 \} $ and $ ( a , b , c , d ) \in \{ 0 , \dots , n \}^4 $, then $ ( m a - n b - 1 )^2 + ( a - c - d - 1 )^2 < 2^{n + 5} $. \end{lemma}

{\bf Proof} The expressions $ m x - n y - 1 $ and $ x - u - v - 1 $ are linear polynomials in $ \mathbb{Z} [ x , y , u , v ] $, so their extrema on $ \{ 0 , \dots , n \}^4 $ can be found by evaluating them at the points $ ( a , b , c , d ) \in \{ 0 , n \}^4 $ only.

By doing so, we find that $ ( m x - n y - 1 )^2 $ (resp., $ ( x - u - v - 1 )^2 $) reach its maximum value $ ( n^2 + 1 )^2 $ (resp., $ ( 2 n + 1 )^2 $) when $ ( x , y ) $ (resp., $ ( x , u , v ) $) is equal to $ ( 0 , n ) $ (resp., $ ( 0 , n , n ) $).

Therefore, the maximum of $ ( m x - n y - 1 )^2 + ( x - u - v - 1 )^2 $ on $ \{ 0 , \dots , n \}^4 $ is $ n^4 + 6 n^2 + 4 n + 2 $, which is strictly upper-bounded by $ 2^{n + 5} $. \qed

\begin{theo}\label{ThmInv} The function $ \inv ( m , n ) $ (for coprime integer arguments $ n \geq 2 $ and $ m \in \{ 1 , \dots , n - 1 \} $) can be represented by the arithmetic term $$ \HW(M(m, n))/(n+5) - (n+1)^4 , $$ where $ M ( m , n ) $ is equal to $$ \begin{array}{llll} \mathcal{C} ( 2 , 4 ) & + \mathcal{A} ( ( m^2 + 1 ) x_1^2 , 4 ) & + \mathcal{A} ( - 2 m n x_1 x_2 , 4 ) & + \mathcal{A} ( - 2 x_1 x_3 , 4 ) \\ & + \mathcal{A} ( - 2 ( m + 1 ) x_1 , 4 ) & + \mathcal{A} ( 2 x_3 x_4 , 4 ) & + \mathcal{A} ( - 2 x_1 x_4 , 4 ) \\ & + \mathcal{A} ( n^2 x_2^2 , 4 ) & + \mathcal{A} ( 2 x_3 , 4 ) & + \mathcal{A} ( x_3^2 , 4 ) \\ & + \mathcal{A} ( 2 n x_2 , 4 ) & + \mathcal{A} ( 2 x_4 , 4 ) & + \mathcal{A} ( x_4^2 , 4 ) . \end{array} $$ \end{theo}

{\bf Proof} It follows from Lemma \ref{LemmaSetInv} and Lemma \ref{LemmaMajInv} by emulating the proof of Lemma \ref{LemmaMazMar}. \qed

\section{The integer part of the root}\label{SectionRoot} 

\begin{lemma}\label{LemmaSetRoot} Given integers $ n \geq 1 $ and $ m \geq 2 $, the number $ \left\lfloor \sqrt[m]{n} \right\rfloor + 1 $ is equal to the cardinality of the set $$ \{ ( a , b ) \in \{ 0 , \dots , n \}^2 : a + b^m - n = 0 \} . $$ \end{lemma}

{\bf Proof} Let $ r = \left\lfloor \sqrt[m]{n} \right\rfloor $, that is to say, the largest integer $ r \geq 0 $ such that $ r^m \leq n < ( r + 1 )^m $.

There are $ r + 1 $ $ m $-th powers that are less or equal to $ n $: $ 0^m $, $ \dots $, $ r^m $.

Hence $ r + 1 $ is the cardinality of the set $ \{ b \in \{ 0 , \dots , n \} : b^m \leq n \} $, which clearly coincides with the cardinality of the set of the statement. \qed

\begin{lemma}\label{LemmaMajRoot1} If $ n $ is a positive integer, $ m $ is an integer exceeding one and $ ( a , b ) \in \{ 0 , \dots , n \}^2 $, then $$ ( a + b^m - n )^2 < 2^{2 m n} . $$ \end{lemma}

{\bf Proof} The largest number of the form $ | a + b^m - n | $, where $ ( a , b ) \in \{ 0 , \dots , n \}^2 $, is clearly $ n^m $. And we have that $ ( n^m )^2 \leq 2^{2 m \lceil \log_2 ( n ) \rceil} < 2^{2 m n} $. \qed

\begin{theo}\label{ThmRootSchematic} Given an integer $ m \geq 2 $, the function $ \left\lfloor \sqrt[m]{n} \right\rfloor $ (for integer arguments $ n \geq 1 $) can be represented by the arithmetic term $$ \HW ( M_m ( n ) ) / ( 2 m n ) - ( n + 1 )^2 - 1 , $$ where $ M_m ( n ) $ is equal to \begin{equation}\label{ExprRoot1} \begin{array}{lll} \mathcal{C} ( n^2 , 2 ) & + \mathcal{A} ( x_1^2 , 2 ) & + \mathcal{A} ( - 2 n x_1 , 2 ) \\ + \mathcal{A} ( 2 x_1 x_2^m , 2 ) & + \mathcal{A} ( x_2^{2 m} , 2 ) & + \mathcal{A} ( - 2 n x_2^m , 2 ) . \end{array} \end{equation} \end{theo}

{\bf Proof} It follows from Lemma \ref{LemmaSetRoot} and Lemma \ref{LemmaMajRoot1} by emulating the proof of Lemma \ref{LemmaMazMar}. \qed

Observe that, as Expression \ref{ExprRoot1} involves the monomial $ x_2^4 $, the arithmetic term from Theorem \ref{ThmRootSchematic} contains the subterm $ G_4 $, which is really cumbersome. Nevertheless, if we are interested just in an arithmetic term representing $ \left\lfloor \sqrt{n} \right\rfloor $, then the proof of Theorem \ref{ThmSquareRoot} will show that we can avoid $ G_4 $ by applying the technique of elimination of arithmetic terms $ G_r $ with $ r \geq 3 $ developed in the proof of Lemma \ref{LemmaExp012}.

\begin{lemma}\label{LemmaMajRoot2} If $ n $ is a positive integer and $ ( a , b , c , d ) \in \{ 0 , \dots , n \}^4 $, then $$ ( a + d - n )^2 + ( b - c )^2 + ( c b - d )^2 < 2^{n + 5} . $$ \end{lemma}

{\bf Proof} Indeed, $$ ( a + d - n )^2 + ( b - c )^2 + ( c b - d )^2 < $$ $$ ( 2 n )^2 + n^2 + ( n^2 )^2 = n^4 + 5 n^2 < 2^{n + 5} . $$ \qed

\begin{theo}\label{ThmSquareRoot} The function $ \left\lfloor \sqrt{n} \right\rfloor $ (for positive integer arguments $ n $) can be represented by the arithmetic term $$ \HW(M(n))/(n+5) - (n+1)^4 - 1 , $$ where $ M ( n ) $ is equal to $$ \begin{array}{llll} \mathcal{C} ( n^2 , 4 ) & + \mathcal{A} ( x_1^2 , 4 ) & + \mathcal{A} ( 2 x_1 x_4 , 4 ) & + \mathcal{A} ( - 2 x_2 x_3 x_4 , 4 ) \\ + \mathcal{A} ( - 2 n x_1 , 4 ) & + \mathcal{A} ( x_2^2 , 4 ) & + \mathcal{A} ( x_2^2 x_3^2 , 4 ) & + \mathcal{A} ( - 2 x_2 x_3 , 4 ) \\ + \mathcal{A} ( - 2 n x_4 , 4 ) & + \mathcal{A} ( x_3^2 , 4 ) & + \mathcal{A} ( 2 x_4^2 , 4 ) . & \end{array} $$ \end{theo}

{\bf Proof} Recall that in the proof of Theorem \ref{ThmRootSchematic} we considered the polynomial $ ( x_1 + x_2^2 - n )^2 $.

The exponent of the only non-exponential occurrence in $ ( x_1 + x_2^2 - n )^2 $ of the variable $ x_2 $ is two, so we consider two new variables $ x_3 $ and $ x_4 $ and the two equations $ x_2 - x_3 = 0 $ and $ x_3 x_2 - x_4 = 0 $, as the proof of Lemma \ref{LemmaExp012} indicates.

Now, we replace in $ ( x_1 + x_2^2 - n )^2 $ the only non-exponential occurrence of $ x_2^2 $ with the variable $ x_4 $, and then we add the squares of the polynomials $ x_2 - x_3 $ and $ x_3 x_2 - x_4 $. The result is $$ ( x_1 + x_4 - n )^2 + ( x_2 - x_3 )^2 + ( x_3 x_2 - x_4 )^2 . $$

The conclusion follows then from Lemma \ref{LemmaSetRoot} and Lemma \ref{LemmaMajRoot2} by emulating the proof of Lemma \ref{LemmaMazMar}. \qed

Note that Theorem \ref{ThmRootSchematic} describes a whole schema of arithmetic terms. It is necessary to fix the value $ m $ beforehand because otherwise the expression $$ n^2 + 2 x_1 x_2^m + x_1^2 - 2 n x_1 + x_2^{2 m} - 2 n x_2^m $$ does not satisfy the definition of simple-in-$ ( x_1 , x_2 ) $ exponential polynomial. The corresponding arithmetic term involves the arithmetic terms $G_{2 m}$ and $ G_m $, which have a different form depending on each chosen $ m $. In addition, the method that was used in Theorem \ref{ThmSquareRoot} for eliminating the arithmetic term $G_4$ cannot be applied in this case, because the number of new variables that are necessary to introduce in order to obtain a new exponential Diophantine definition is not independent of $ m $. However, Theorem \ref{ThmRootUniform} will provide a single simple-in-$ ( x_1 , \dots , x_7 ) $ exponential polynomial in which $ m $ occurs as a variable.

In what remains of the present section, we denote the expression $$ ( x_1 - ( m + 1 ) x_7 - 1 )^2 + ( x_2 - m x_1 )^2 + ( 2^{x_2} - x_3 2^{x_1} + x_3 x_7 - x_4 )^2 + ( x_4 + x_5 - 2^{x_1} + x_7 + 1 )^2 + ( x_4 + x_6 - n )^2 $$ by $ E ( m , n , x_1 , \dots , x_7 ) $.

\begin{lemma}\label{LemmaSetRootUniform} Given two integers $ m \geq 2 $ and $ n \geq 1 $, the number $ \lfloor \sqrt[m]{n} \rfloor + 1 $ is equal to the cardinality of the set $$ \{ ( a , b , c , d , e , f , g ) \in \{ 0 , \dots , 2^{nm^2 + nm + 1} - 1 \}^7 : E ( m , n , a , b , c , d , e , f , g ) = 0 \} . $$ \end{lemma}

{\bf Proof} Let $ A $ be the set of the statement and let $ B $ be the set of points of the form $$ \left( g m + g + 1 , ( g m + g + 1 ) m , \left \lfloor \frac{2^{(gm + g + 1)m}}{2^{gm + g + 1} - g} \right \rfloor , \right. $$ $$ 2^{(mg + g + 1)m} - (2^{mg + g + 1} - g) \left \lfloor \frac{2^{(gm + g + 1)m}}{2^{gm + g + 1} - g} \right \rfloor , $$ $$ - 2^{(mg + g + 1)m} + (2^{mg + g + 1} - g) \left \lfloor \frac{2^{(gm + g + 1)m}}{2^{gm + g + 1} - g} \right \rfloor + 2^{mg + g + 1} - g - 1 , $$ $$ \left. n - 2^{(mg + g + 1)m} + (2^{mg + g + 1} - g) \left \lfloor \frac{2^{(gm + g + 1)m}}{2^{gm + g + 1} - g} \right \rfloor , g \right) , $$ where $ g \in \{ 0 , \dots , 2^{nm^2 + nm + 1} - 1 \} $. 

The inclusion $ B \subseteq A $ is immediate.

The cardinality of $ B $ is $ \left \lfloor \sqrt[m]{n}\right \rfloor + 1 $: indeed, notice that, because of Identity \ref{EqMarchenkov}, the following equalities hold: $$ 2^{(mg + g + 1)m} - (2^{mg + g + 1} - g) \left \lfloor \frac{2^{(gm + g + 1)m}}{2^{gm + g + 1} - g} \right \rfloor = 2^{( g m + g + 1 ) m} \bmod ( 2^{g m + g + 1} - g ) = g^m . $$ Hence, by applying the condition $ x_4 + x_6 - n = 0 $ (which must be satisfied because $ B \subseteq A $), we have that $ g^m \leq n $. And the non-negative integers $ g $ such that $ g^m \leq n $ are exactly the elements of the set $\{ 0 , 1 , \dots , \lfloor \sqrt[m]{n} \rfloor \} $, whose cardinality is $ \left \lfloor \sqrt[m]{n}\right \rfloor + 1 $.

It only remains to show that $ A \subseteq B $, so let $ ( a , b , c , d , e , f , g ) \in A $.

The conditions $ a - ( m + 1 ) g - 1 = 0 $ and $ b - m a = 0 $ imply that $ a = g m + g + 1 $ and $ b = ( g m + g + 1 ) m $.

The condition $ 2^b - c 2^a + c g - d = 0 $ implies that $ 2^b = c ( 2^a - g ) + d $. And the condition $ d + e - 2^a + g + 1 = 0 $ implies that $ d < 2^a - g $, from which follows that $$ c = \left \lfloor \frac{2^b}{2^a - g} \right \rfloor = \left \lfloor \frac{2^{(gm + g + 1)m}}{2^{gm + g + 1} - g} \right \rfloor . $$

Thus $$ d = 2^b - c 2^a + c g = 2^{(mg + g + 1)m} - (2^{mg + g + 1} - g) \left \lfloor \frac{2^{(gm + g + 1)m}}{2^{gm + g + 1} - g} \right \rfloor $$ and, consequently, $$ e = - d + 2^a - g - 1 = - 2^{(mg + g + 1)m} + (2^{mg + g + 1} - g) \left \lfloor \frac{2^{(gm + g + 1)m}}{2^{gm + g + 1} - g} \right \rfloor + 2^{mg + g + 1} - g - 1 . $$

Finally, the condition $ d + f - n = 0 $ implies that $$ f = n - 2^{(mg + g + 1)m} + (2^{mg + g + 1} - g) \left \lfloor \frac{2^{(gm + g + 1)m}}{2^{gm + g + 1} - g} \right \rfloor . $$ \qed

\begin{lemma}\label{LemmaMajRootUniform} Given two integers $ m \geq 2 $ and $ n \geq 1 $, and a point $$ ( a , b , c , d , e , f , g ) \in \{ 0 , \dots , 2^{nm^2 + nm + 1} - 1 \}^7 , $$ we have that $$ E ( m , n , a , b , c , d , e , f , g ) < 2^{2^{nm^2 + nm + 2} + 2 (nm^2 + nm) + 9}. $$ \end{lemma}

{\bf Proof} Let $ F ( m , n , x_1 , \dots , x_7 ) $ be the expression that is obtained by replacing every minus sign with a plus sign in the expression $ E ( m , n , x_1 , \dots , x_7 ) $, and let $ t = 2^{nm^2 + nm + 1} $.

Then $$ F ( t , \dots , t ) = ( t + ( t + 1 ) t + 1 )^2 + ( t + t^2 )^2 + ( 2^t + t 2^t + t^2 + t )^2 + ( 3 t + 2^t + 1 )^2 + 9 t^2 < $$ $$ 5 ( 2^t + t 2^t + t^2 + t )^2 < 5 ( 4 t 2^t )^2 < 2^7 t^2 2^{2 t} = 2^{2^{nm^2 + nm + 2} + 2 (nm^2 + nm) + 9} . $$ \qed

\begin{theo}\label{ThmRootUniform} The function $ \left\lfloor \sqrt[m]{n} \right\rfloor $ (for integer arguments $ n \geq 1 $ and $ m \geq 2 $) can be represented by the arithmetic term $$ \HW ( M ( m , n ) ) / ( 2^{nm^2 + nm + 2} + 2 (nm^2 + nm) + 9 ) - 2^{7 (nm^2 + nm + 1) } - 1 , $$ where $ M ( m , n ) $ is equal to $$ \begin{array}{llll} \mathcal{C} ( 2 + n^2 , 7 ) & + \mathcal{A} ( - 2 m x_1 x_2 , 7 ) & + \mathcal{A} ( ( m^2 + 1 ) x_1^2 , 7 ) & + \mathcal{A} ( - x_4 2^{x_1 + 1} , 7 ) \\ + \mathcal{A} ( - 2 x_1 , 7 ) & + \mathcal{A} ( - 2 ( m + 1 ) x_1 x_7 , 7 ) & + \mathcal{A} ( x_2^2 , 7 ) & + \mathcal{A} ( - x_5 2^{x_1 + 1} , 7 ) \\ + \mathcal{A} ( 2 x_5 , 7 ) & + \mathcal{A} ( 2 x_4 x_5 , 7 ) & + \mathcal{A} ( 3 x_4^2 , 7 ) & + \mathcal{A} ( - x_7 2^{x_1 + 1} , 7 ) \\ + \mathcal{A} ( 2 ( 1 - n ) x_4 , 7 ) & + \mathcal{A} ( 2 x_4 x_6 , 7 ) & + \mathcal{A} ( x_5^2 , 7 ) & + \mathcal{A} ( - x_4 2^{x_2 + 1} , 7 ) \\ + \mathcal{A} ( - 2 n x_6 , 7 ) & + \mathcal{A} ( 2 x_4 x_7 , 7 ) & + \mathcal{A} ( x_6^2 , 7 ) & + \mathcal{A} ( x_3 x_4 2^{x_1 + 1} , 7 ) \\ + \mathcal{A} ( 2 ( m + 2 ) x_7 , 7 ) & + \mathcal{A} ( 2 x_5 x_7 , 7 ) & + \mathcal{A} ( ( m^2 + 2 m + 2 ) x_7^2 , 7 ) & + \mathcal{A} ( x_3 x_7 2^{x_2 + 1} , 7 ) \\ + \mathcal{A} ( - 2 x_3 x_4 x_7 , 7 ) & + \mathcal{A} ( - 2^{x_1 + 1} , 7 ) & + \mathcal{A} ( 2^{2 x_1} , 7 ) & + \mathcal{A} ( - x_3^2 x_7 2^{x_1 + 1} , 7 ) \\ + \mathcal{A} ( x_3^2 x_7^2 , 7 ) & + \mathcal{A} ( x_3^2 2^{2 x_1} , 7 ) & + \mathcal{A} ( 2^{2 x_2} , 7 ) & + \mathcal{A} ( - x_3 2^{x_1 + x_2 + 1} , 7 ) . \end{array} $$ \end{theo}

{\bf Proof} It follows from Lemma \ref{LemmaSetRootUniform} and Lemma \ref{LemmaMajRootUniform} by emulating the proof of Lemma \ref{LemmaMazMar}. \qed

Finally, we show two applications of Theorem \ref{ThmSquareRoot}.

A \textbf{semiprime} is a product of exactly two primes (see Weisstein \cite{Weisstein2}). By combining Theorem \ref{ThmPhi} with Theorem \ref{ThmSquareRoot} we get Theorem \ref{ThmFactorizationRSA}, an unexpected byproduct on squarefree semiprimes, which are the moduli for the public encryption method RSA (cf.\ Weisstein \cite{Weisstein1}).

\begin{theo}\label{ThmFactorizationRSA} 
There is an arithmetic term $ T ( n ) $ such that $ T ( p q ) = q $ for every two primes $ p $ and $ q $ such that $ p < q $. 
\end{theo}

{\bf Proof} Let $ N = p q $.

The following trick is folklore in public-key cryptography: $ \varphi(N) $ is equal to $ (p-1)(q-1) $ (recall Identity \ref{EqCharPhi}) or, equivalently, to $ N - (q+p) + 1 $; so \begin{equation}\label{EqSum} q + p = N - \varphi(N) + 1 . \end{equation}

In addition, $ (q-p)^2 $ is equal to $ (q+p)^2 - 4pq $ or, in other words, to $ (N - \varphi(N) + 1)^2 - 4N $; from which follows that \begin{equation}\label{EqSubtr} q - p = \lfloor \sqrt{(N - \varphi(N) + 1)^2 - 4N} \rfloor . \end{equation}

Therefore, by summing Identity \ref{EqSum} and Identity \ref{EqSubtr}, we conclude that $$ q = ( N - \varphi(N) + 1 + \lfloor \sqrt{(N - \varphi(N) + 1)^2 - 4N} \rfloor ) / 2 . $$ \qed

The \textbf{factoring problem} consist in, given any integer $ n > 1 $, finding an integer $ d > 1 $ such that $ d $ divides $ n $ (cf.\ Nederlof \cite[Example 2]{Nederlof}).

Note that, as the function that associate each integer $ n > 1 $ to its least prime divisor is also a Kalmar function, there exists an arithmetic term which outputs the least prime factor of $ n $. Such an arithmetic term remains to be found, but its existence solves (although probably not efficiently) the factoring problem.

Another application of Theorem \ref{ThmSquareRoot} is Theorem \ref{ThmCantorsPairing}, which provides an arithmetic term that represents \textbf{Cantor's pairing function}, the bijection that maps each pair $ ( x , y )$ of non-negative integers into the non-negative integer $ (x+y)(x+y+1) / 2 + x $ (see Weisstein \cite{Weisstein0}).

\begin{theo}\label{ThmCantorsPairing} If $ c $ is Cantor's pairing function, then there are arithmetic terms $x(n)$ and $y(n)$ such that $ c(x(n), y(n)) = n $. \end{theo}

{\bf Proof} The usual computation method is the following (cf.\ Weisstein \cite{Weisstein0}): \begin{eqnarray*} w(n) &:=& \lfloor ( \sqrt{8n+1} - 1 ) / 2 \rfloor , \\ t(n) &:=& ( w(n)^2 + w(n) ) / 2 , \\ x(n) &:=& n - t(n) , \\ y(n) &:=& w(n) - x(n) . \end{eqnarray*}

So it only remains to write $w(n)$ as an arithmetic term.

The identity $ w(n) = \lfloor ( \sqrt{8n+1} - 1 ) / 2 \rfloor $ is equivalent with $ 2w(n) + 1 \leq \sqrt{8n+1} < 2 w(n) + 3 $.

If $\lfloor \sqrt{8 n + 1} \rfloor$ is odd (resp., even), then it is equal to $ 2w(n) + 1$ (resp., $ 2w(n) + 2 $) and consequently $ w(n) $ equals $ ( \lfloor \sqrt{8 n + 1} \rfloor - 1 ) / 2 $ (resp., $ ( \lfloor \sqrt{8 n + 1} \rfloor - 2 ) / 2 $).

Therefore, $w(n)$ is represented by the arithmetic term $ ( \lfloor \sqrt{8 n + 1} \rfloor - 2 + ( \lfloor \sqrt{8 n + 1} \rfloor \bmod 2 ) ) / 2 $. \qed

\section{The integer part of the logarithm}\label{SectionLog}

\begin{lemma}\label{LemmaSetLog} Given two integers $ m \geq 2 $ and $ n \geq 1 $, the number $ \left\lfloor \log_m ( n ) \right\rfloor + 1 $ is equal to the cardinality of the set $$ \{ ( a , b ) \in \{ 0 , \dots , n \}^2 : a + m^b - n = 0 \} . $$ \end{lemma}

{\bf Proof} Let $ r = \left\lfloor \log_m ( n ) \right\rfloor $, that is to say, the only element $ r \in \{ 0 , \dots , n \} $ such that $ m^r \leq n < m^{r + 1} $.

There are $ r + 1 $ powers of $ m $ that do not exceed $ n $: $ m^0 $, $m ^1 $, $ \dots $, $ m^r $.

Hence $ r + 1 $ is the cardinality of the set $ \{ b \in \{ 0 , \dots , n \} : m^b \leq n \} $, which clearly coincides with the cardinality of the set of the statement. \qed

\begin{lemma}\label{LemmaMajLog} If $ m $ is an integer exceeding one, $ n $ is a positive integer and $ ( a , b ) \in \{ 0 , \dots , n \}^2 $, then $ ( a + m^b - n )^2 < 2^{2 m n} $. \end{lemma}

{\bf Proof} The largest number of the form $ | a + m^b - n | $, where $ ( a , b ) \in \{ 0 , \dots , n \}^2 $, is clearly $ m^n $.

And we have that $ ( m^n )^2 \leq 2^{2 n \lceil \log_2 ( m ) \rceil} < 2^{2 m n} $. \qed

\begin{theo}\label{ThmLog} The function $ \left\lfloor \log_m ( n ) \right\rfloor $ (for integer arguments $ m \geq 2 $ and $ n \geq 1 $) can be represented by the arithmetic term $$ \HW(M( m , n ))/( 2 m n ) - ( n + 1 )^2 - 1 , $$ where $ M ( m , n ) $ is equal to $$ \begin{array}{ll} \mathcal{C} ( n^2 , 2 ) & + \mathcal{A} ( 2 x_1 m^{x_2} , 2 ) \\ + \mathcal{A} ( - 2 n m^{x_2} , 2 ) & + \mathcal{A} ( m^{2 x_2} , 2 ) \\ + \mathcal{A} ( - 2 n x_1 , 2 ) & + \mathcal{A} ( x_1^2 , 2 ) . \end{array} $$ \end{theo}

{\bf Proof} It follows from Lemma \ref{LemmaSetLog} and Lemma \ref{LemmaMajLog} by emulating the proof of Lemma \ref{LemmaMazMar}. \qed

Notice that the exponent $ n + 1 $ used in Theorem \ref{ThmPadicVal} is, in general, far too big. Theorem \ref{ThmEfficientPadicVal} provides a more efficient arithmetic term for the $ p $-adic valuation.

\begin{theo}\label{ThmEfficientPadicVal} The function $ \nu_p ( n ) $ (for integer arguments $ n \geq 1 $ and $ p $ prime) can be represented by the arithmetic term $$ \left\lfloor \dfrac{\gcd \left( n , p^{\lfloor \log_p ( n ) \rfloor + 1} \right)^{\lfloor \log_p ( n ) \rfloor + 3} \bmod \left( p^{\lfloor \log_p ( n ) \rfloor + 3} - 1 \right)^2}{p^{\lfloor \log_p ( n ) \rfloor + 3} - 1} \right\rfloor . $$ \end{theo}

{\bf Proof} Let $ x = \nu_p ( n ) $.

It is clear that $$ x < n \leq p^{\lceil \log_p ( n ) \rceil} < p^{\lceil \log_p ( n ) \rceil + 1} = p^{\lfloor \log_p ( n ) \rfloor + 2} , $$ so $ x < p^{\lfloor \log_p ( n ) \rfloor + 3} - 1 $ and consequently $ 1 + x \left( p^{\lfloor \log_p ( n ) \rfloor + 3} - 1 \right) < \left( p^{\lfloor \log_p ( n ) \rfloor + 3} - 1 \right)^2 $.

In addition, $$ \gcd \left( n , p^{\lceil \log_p ( n ) \rceil} \right)^{\lfloor \log_p ( n ) \rfloor + 3} = {\left ( p^x \right )}^{\lfloor \log_p ( n ) \rfloor + 3} = {\left ( p^{\lfloor \log_p ( n ) \rfloor + 3} \right )}^{x} = $$ $$ {\left ( p^{\lfloor \log_p ( n ) \rfloor + 3} - 1 + 1 \right )}^{x} = \sum_{k=0}^x \left( \binom{x}{k} \left( p^{\lfloor \log_p ( n ) \rfloor + 3} - 1 \right)^k \right) = $$ $$ 1 + x \left( p^{\lfloor \log_p ( n ) \rfloor + 3} - 1 \right) + \left( p^{\lfloor \log_p ( n ) \rfloor + 3} - 1 \right)^2 \sum_{k = 2}^x \left( \binom{x}{k} \left( p^{\lfloor \log_p ( n ) \rfloor + 3} - 1 \right)^{k-2} \right) . $$

Thus $$ \gcd \left( n , p^{\lceil \log_p ( n ) \rceil} \right)^{\lfloor \log_p ( n ) \rfloor + 3} \bmod \left( p^{\lfloor \log_p ( n ) \rfloor + 3} - 1 \right)^2 = 1 + x \left( p^{\lfloor \log_p ( n ) \rfloor + 3} - 1 \right) , $$ from which the statement immediately follows. \qed

\section{The multiplicative order}\label{SectionOrd}

Theorem \ref{EulersTheorem} is known as \textbf{Euler's theorem} (see Rosen \cite[Theorem 6.14]{Rosen}).

\begin{theo}\label{EulersTheorem} Given two coprime integers $ n \geq 2 $ and $ m $, then $ m^{\varphi ( n )} \equiv 1 \pmod{n} $. \end{theo}

Lemma \ref{LemmaOrd} is an instance of Rosen \cite[Corollary 9.1.1]{Rosen}.

\begin{lemma}\label{LemmaOrd} Given two coprime integers $ n \geq 2 $ and $ m \in \{ 1 , \dots , n - 1 \} $, we have that $ \ord ( m , n ) $ divides $ \varphi ( n ) $. \end{lemma}

Lemma \ref{LemmaFactorization} is also of frequent use (see, for example, Sauras-Altuzarra \cite[Lemma 3.3.3.7]{SaurasAltuzarra}).

\begin{lemma}\label{LemmaFactorization} Given two coprime integers $ n \geq 2 $ and $ m \in \{ 1 , \dots , n - 1 \} $, and one integer $ r > 0 $ such that $ n $ divides $ m^r - 1 $, we have that $ \ord ( m , n ) $ divides $ r $. \end{lemma}

{\bf Proof} Suppose the contrary.

Then there are two positive integers $ x $ and $ y $ such that $ x \ord ( m , n ) + y = r $ and $ y < \ord ( m , n ) $ (because, by definition, $ \ord ( m , n ) \leq r $).

Therefore $$ 1 \equiv m^r = m^{x \ord ( m , n ) + y} = ( m^{\ord ( m , n )} )^x m^y \equiv 1^x m^y = m^y \ \textrm{(mod $ n $)} $$ (because $ n $ divides $ m^r - 1 $), which contradicts the fact that $ \ord ( m , n ) $ is the minimum positive integer $ k $ such that $ n $ divides $ m^k - 1 $. \qed

\begin{lemma}\label{LemmaSetOrd} Given two coprime integers $ n \geq 2 $ and $ m \in \{ 1 , \dots , n - 1 \} $, the number $$ \dfrac{\varphi ( n )}{\ord ( m , n )} $$ is equal to the cardinality of the set $$ \{ ( a , b , c , d ) \in \{ 0 , \dots , m^{\varphi ( n )} \}^4 : ( m^a - n b - 1 )^2 + ( a - c - 1 )^2 + ( \varphi ( n ) - d - a )^2 = 0 \} . $$ \end{lemma}

{\bf Proof} Let $ A $ be the set of the statement, let $ x = \ord ( m , n ) $ and let $ B $ the set $$ \{ ( u x , ( m^{u x} - 1 ) / n , u x - 1 , \varphi ( n ) - u x ) : u \in \{ 1 , \dots , \varphi ( n ) / x \} \} . $$

The cardinality of $ B $ is obviously $ \varphi ( n ) / x $.

We know that $ m^x \equiv 1 \pmod{n} $, so $ m^{ux} \equiv 1 \pmod{n} $ and thus $ ( m^{ux} - 1 ) / n $ is an integer.

Therefore $ B \subseteq \{ 0 , \dots , m^{\varphi ( n )} \}^4 $, so clearly $ B \subseteq A $.

It only remains to check that $ A \subseteq B $. Let $ ( a , b , c , d ) \in A $.

The condition $ m^a - n b - 1 = 0 $ implies that $ n $ divides $ m^a - 1 $, so $ x | a $ by applying Lemma \ref{LemmaFactorization}.

Hence there is some number $ u \in \{ 1 , \dots , m^{\varphi ( n )} \} $ such that $ a = u x $.

The condition $ \varphi(n) - d - a = 0 $ implies that $ a \leq \varphi(n) $, so $ u \in \{ 1 , \dots , \varphi ( n ) / x \} $. 

By again applying the condition $ m^a - n b - 1 = 0 $, we get that $ b = ( m^{u x} - 1 ) / n $.

And finally, the conditions $ a - c - 1 = 0 $ and $ \varphi ( n ) - d - a = 0 $ imply that $ c = u x - 1 $ and $ d = \varphi ( n ) - u x $. \qed

\begin{lemma}\label{LemmaMajOrd} Given two coprime integers $ n \geq 3 $ and $ m \in \{ 2 , \dots , n - 1 \} $, and a point $$ ( a , b , c , d ) \in \{ 0 , \dots , m^{\varphi ( n )} \}^4 , $$ we have that $$ ( m^a - n b - 1 )^2 + ( a - c - 1 )^2 + ( \varphi ( n ) - d - a )^2 < 2^{2 m^{\varphi ( n ) + 1} + 2} . $$ \end{lemma}

{\bf Proof} Indeed, $$ ( m^a - n b - 1 )^2 + ( a - c - 1 )^2 + ( \varphi ( n ) - d - a )^2 \leq $$ $$ ( m^{m^{\varphi ( n )}} - 1 )^2 + ( m^{\varphi ( n )} + 1 )^2 + ( 2 m^{\varphi ( n )} - \varphi ( n ) )^2 < $$ $$ 3 ( m^{m^{\varphi ( n )}} )^2 < ( 2 \cdot 2^{m^{\varphi ( n )} \lceil \log_2 ( m ) \rceil} )^2 < ( 2^{1 + m^{\varphi ( n ) + 1}} )^2 < 2^{2 m^{\varphi ( n ) + 1} + 2} . $$ \qed

Theorem \ref{ThmOrd} provides an arithmetic term for the function $ \ord $ in terms of the function $ \varphi $, whose arithmetic-term representation has been already shown in Theorem \ref{ThmPhi}.

\begin{theo}\label{ThmOrd} The function $ \ord ( m , n ) $ (for coprime integer arguments $ n \geq 3 $ and $ m \in \{ 2 , \dots , n - 1 \} $) can be represented by the arithmetic term $$ \dfrac{\varphi ( n )}{\HW(M(m, n)) / ( 2 m^{\varphi ( n ) + 1} + 2 ) - ( m^{\varphi(n)} + 1 )^4} , $$ where $ M ( m , n ) $ is equal to $$ \begin{array}{llll} \mathcal{C} ( \varphi ( n )^2 + 2 , 4 ) & + \mathcal{A} ( 2 n x_2 , 4 ) & + \mathcal{A} ( 2 x_1^2 , 4 ) & + \mathcal{A} ( - 2 ( \varphi ( n ) + 1 ) x_1 , 4 ) \\ + \mathcal{A} ( - 2 n m^{x_1} x_2 , 4 ) & + \mathcal{A} ( 2 x_3 , 4 ) & + \mathcal{A} ( n^2 x_2^2 , 4 ) & + \mathcal{A} ( - 2 \varphi ( n ) x_4 , 4 ) \\ + \mathcal{A} ( - 2 m^{x_1} , 4 ) & + \mathcal{A} ( - 2 x_1 x_3 , 4 ) & + \mathcal{A} ( x_3^2 , 4 ) & \\ + \mathcal{A} ( m^{2 x_1} , 4 ) & + \mathcal{A} ( 2 x_1 x_4 , 4 ) & + \mathcal{A} ( x_4^2 , 4 ) . & \end{array} $$ \end{theo}

{\bf Proof} It follows from Lemma \ref{LemmaSetOrd} and Lemma \ref{LemmaMajOrd} by emulating the proof of Lemma \ref{LemmaMazMar}. \qed

\section{The discrete logarithm}\label{SectionDiscLog}

\begin{lemma}\label{LemmaSetDiscLog} Given two coprime integers $ n \geq 3 $ and $ m \in \{ 2 , \dots , n - 1 \} $, and a primitive root $ g $ modulo $ n $ such that $ g \geq 2 $, the number $ \dlog ( m , g , n ) $ is equal to the cardinality of the set $$ \{ ( a , b , c , d ) \in \{ 0 , \dots , g^{\varphi ( n )} \}^4 : ( a + b + c + 1 - \varphi ( n ) )^2 + ( g^{a + b + 1} - n d - m )^2 = 0 \} . $$ \end{lemma}

{\bf Proof} Let $ x $ be the number $ \dlog ( m , g , n ) $, let $ A $ be the set of the statement and let $ B $ be the set $$ \{ ( a , x - a - 1 , \varphi ( n ) - x , ( g^x - m ) / n ) : a \in \{ 0 , \dots , g^{\varphi ( n )} \} \ \textrm{and} \ x - a - 1 \geq 0 \} . $$

The number $ x - a - 1 $ is non-negative, so $ a \in \{ 0 , \dots , x - 1 \} $ and hence the cardinality of $ B $ is $ x $.

The inclusion $ B \subseteq A $ is immediate, so it only remains to check that $ A \subseteq B $. Let $ ( a , b , c , d ) \in A $.

From the condition $ a + b + c + 1 - \varphi ( n ) = 0 $ we get that $ a + b + 1 \leq \varphi ( n ) $.

And, from the condition $ g^{a + b + 1} - n d - m = 0 $ we get that $ g^{a + b + 1} \equiv m \pmod{n} $.

It follows that $ x = a + b + 1 = \varphi ( n ) - c $, so $ b = x - a - 1 $, $ c = \varphi ( n ) - x $ and $ d = ( g^x - m ) / n $. \qed

\begin{lemma}\label{LemmaMajDiscLog} Given two coprime integers $ n \geq 3 $ and $ m \in \{ 2 , \dots , n - 1 \} $, a primitive root $ g $ modulo $ n $ such that $ g \geq 2 $ and a point $ ( a , b , c , d ) \in \{ 0 , \dots , g^{\varphi ( n )} \}^4 $, we have that $$ ( a + b + c + 1 - \varphi ( n ) )^2 + ( g^{a + b + 1} - n d - m )^2 < 2^{5 + 2 n + 2 g ( 2 g^{\varphi ( n )} + 1 )} . $$ \end{lemma}

{\bf Proof} Indeed, $$ ( a + b + c + 1 - \varphi ( n ) )^2 + ( g^{a + b + 1} - n d - m )^2 < $$ $$ ( a + b + c + 1 + \varphi ( n ) )^2 + ( g^{a + b + 1} + n d + m )^2 < $$ $$ ( 5 g^{\varphi ( n )} )^2 + ( 3 n g^{2 g^{\varphi ( n )} + 1} )^2 < 2 ( 3 n g^{2 g^{\varphi ( n )} + 1} )^2 < $$ $$ 2^5 n^2 g^{2 ( 2 g^{\varphi ( n )} + 1 )} \leq 2^{5 + 2 \lceil \log_2 ( n ) \rceil + 2 ( 2 g^{\varphi ( n )} + 1 ) \lceil \log_2 ( g ) \rceil} < 2^{5 + 2 n + 2 g ( 2 g^{\varphi ( n )} + 1 )} . $$ \qed

Like Theorem \ref{ThmFactorizationRSA} and Theorem \ref{ThmOrd}, Theorem \ref{ThmDiscLog} makes use of the arithmetic term $ \varphi ( n ) $ from Theorem \ref{ThmPhi}.

\begin{theo}\label{ThmDiscLog} The function $ \dlog ( m , g , n ) $ (for integer arguments $ n \geq 3 $, $ m \in \{ 2 , \dots , n - 1 \} $ and $ g \geq 2 $ such that $ n $ and $ m $ are coprime and $ g $ is a primitive root modulo $ n $) can be represented by the arithmetic term $$ \HW ( M ( m , g , n ) ) / ( 5 + 2 n + 2 g ( 2 g^{\varphi ( n )} + 1 ) ) - ( g^{\varphi ( n )} + 1 )^4 , $$ where $ M ( m , g , n ) $ is equal to $$ \begin{array}{lll} \mathcal{C} ( 1 + m^2 - 2 \varphi ( n ) + \varphi ( n )^2 , 4 ) & + \mathcal{A} ( 2 m n x_4 , 4 ) & + \mathcal{A} ( n^2 x_4^2 , 4 ) \\ + \mathcal{A} ( 2 ( 1 - \varphi ( n ) ) x_1 , 4 ) & + \mathcal{A} ( - 2 m g^{x_1 + x_2 + 1} , 4 ) & + \mathcal{A} ( x_1^2 , 4 ) \\ + \mathcal{A} ( 2 ( 1 - \varphi ( n ) ) x_2 , 4 ) & + \mathcal{A} ( g^{2 x_1 + 2 x_2 + 2} , 4 ) & + \mathcal{A} ( x_2^2 , 4 ) \\ + \mathcal{A} ( 2 ( 1 - \varphi ( n ) ) x_3 , 4 ) & + \mathcal{A} ( - 2 n x_4 g^{x_1 + x_2 + 1} , 4 ) & + \mathcal{A} ( x_3^2 , 4 ) \\ + \mathcal{A} ( 2 x_1 x_2 , 4 ) & + \mathcal{A} ( 2 x_1 x_3 , 4 ) & + \mathcal{A} ( 2 x_2 x_3 , 4 ) . \end{array} $$ \end{theo}

{\bf Proof} It follows from Lemma \ref{LemmaSetDiscLog} and Lemma \ref{LemmaMajDiscLog} by emulating the proof of Lemma \ref{LemmaMazMar}. \qed

\begin{appendices}

\section{Maple codes}\label{MapleCodes}

Most of the verification process is based on the following Maple code, to which we refer as the \textbf{base code}. In it we define the Hamming weight, the generalized geometric progressions and the functions $ \mathcal{C} $ and $ \mathcal{A} $ from the representation method explained in Section \ref{SectionRepresentationMethod}.

\begin{Verbatim}
HW := n -> add(convert(n, base, 2)):
G[0] := (q, t) -> (q^(t+1)-1)/(q-1):
G[1] := (q, t) -> q*(t*q^(t+1)-(t+1)*q^t+1)/(q-1)^2:
G[2] := (q, t) -> q*(t^2*q^(t+2)-(2*t^2+2*t-1)*q^(t+1)+(t+1)^2*q^t-q-1)/(q-1)^3:
C := (e, k, t, w) ->
   (2^w-e+1)*(2^(2*w*t^k)-1)/(2^w+1):
A := (a, U, B, V, k, t, w) ->
   -(2^w-1)*a*mul(G[U[i]](B[i]^V[i]*2^(2*w*t^(i-1)), t-1), i = 1 .. k):
\end{Verbatim}

Therefore, the expressions of the form $$ \mathcal{C} ( e ( \vec{n} ) , k ) , $$ $$ \mathcal{A} ( a ( \vec{n} ) x_1^{u_1} \dots x_k^{u_k} b_1 ( \vec{n} )^{v_1 ( \vec{n} ) x_1} \dots b_k ( \vec{n} )^{v_k ( \vec{n} ) x_k} , k ) $$ are encoded, respectively, as
\begin{Verbatim}
C(e, k, t, w),
\end{Verbatim}
\begin{Verbatim}
A(a, U, B, V, k, t, w).
\end{Verbatim}

Observe that, in the base code, we utilize a quick Maple command in order to define the Hamming weight. For displaying the arithmetic term representing the Hamming weight, which is extremely inefficient, we can use the following code. It defines the greatest common divisor, the dyadic valuation and the Hamming weight by utilizing the arithmetic terms from Section \ref{SectionUsefulTerms}.

\begin{Verbatim}
gcd2 := (m, n) ->
   irem(floor(((2^(m^2*n*(n+1))-2^(m^2*n))*(2^(m^2*n^2)-1))
   /((2^(m^2*n)-1)*(2^(m*n^2)-1)*2^(m^2*n^2))), 2^(m*n)):
nu2 := n -> floor(irem(gcd2(n, 2^n)^(n+1), (2^(n+1)-1)^2)/(2^(n+1)-1)):
HW := n -> nu2(irem(floor((1+2^(2*n))^(2*n)/2^(2n^2)), 2^(2*n))):
lprint(HW(n));
\end{Verbatim}

The previous code produces the following output (of 1039 characters).

\begin{Verbatim}
floor(irem(irem(floor((2^(irem(floor((1+2^(2*n))^(2*n)
/2^(2*n^2)), 2^(2*n))^2*2^irem(floor((1+2^(2*n))^(2*n)
/2^(2*n^2)), 2^(2*n))*(2^irem(floor((1+2^(2*n))^(2*n)
/2^(2*n^2)), 2^(2*n))+1))-2^(irem(floor((1+2^(2*n))^(2*n)
/2^(2*n^2)), 2^(2*n))^2*2^irem(floor((1+2^(2*n))^(2*n)
/2^(2*n^2)), 2^(2*n))))*(2^(irem(floor((1+2^(2*n))^(2*n)
/2^(2*n^2)), 2^(2*n))^2*(2^irem(floor((1+2^(2*n))^(2*n)
/2^(2*n^2)), 2^(2*n)))^2)-1)/((2^(irem(floor((1+2^(2*n))^(2*n)
/2^(2*n^2)), 2^(2*n))^2*2^irem(floor((1+2^(2*n))^(2*n)
/2^(2*n^2)), 2^(2*n)))-1)*(2^(irem(floor((1+2^(2*n))^(2*n)
/2^(2*n^2)), 2^(2*n))*(2^irem(floor((1+2^(2*n))^(2*n)
/2^(2*n^2)), 2^(2*n)))^2)-1)*2^(irem(floor((1+2^(2*n))^(2*n)
/2^(2*n^2)), 2^(2*n))^2*(2^irem(floor((1+2^(2*n))^(2*n)
/2^(2*n^2)), 2^(2*n)))^2))), 2^(irem(floor((1+2^(2*n))^(2*n)
/2^(2*n^2)), 2^(2*n))*2^irem(floor((1+2^(2*n))^(2*n)
/2^(2*n^2)), 2^(2*n))))^(irem(floor((1+2^(2*n))^(2*n)
/2^(2*n^2)), 2^(2*n))+1), (2^(irem(floor((1+2^(2*n))^(2*n)
/2^(2*n^2)), 2^(2*n))+1)-1)^2)/(2^(irem(floor((1+2^(2*n))^(2*n)
/2^(2*n^2)), 2^(2*n))+1)-1))
\end{Verbatim}

We can use the following code in order to test Lemma \ref{LemmaSetTau}.

\begin{Verbatim}
TestTau := proc(n)
   local a, b, L:
   L := []:
   for a from 0 to n do
      for b from 0 to n do
         if n-a*b = 0
         then L := [op(L), [a, b]]: fi: od: od:
   L: end:
seq(numtheory:-tau(n), n = 1 .. 25);
seq(nops(TestTau(n)), n = 1 .. 25);
\end{Verbatim}

We can experimentally verify Theorem \ref{ThmTau} by extending the base code with the following one.

\begin{Verbatim}
k := 2:
t := n -> n+1:
w := n -> n+4:
M := n ->
   C(n^2, k, t(n), w(n))
   + A(-2*n, [1, 1], [2, 2], [0, 0], k, t(n), w(n))
   + A(1, [2, 2], [2, 2], [0, 0], k, t(n), w(n)):
Tau := n -> HW(M(n))/w(n)-t(n)^k:
seq(numtheory:-tau(n), n = 1 .. 25);
seq(Tau(n), n = 1 .. 25);
\end{Verbatim}

We can use the following code in order to test Lemma \ref{LemmaSetSigma}.

\begin{Verbatim}
TestSigma := proc(n)
   local a, b, c, L:
   L := []:
   for a from 0 to n do
      for b from 0 to n do
         for c from 0 to n do
            if n-(a+b+1)*c = 0
            then L := [op(L), [a, b, c]]: fi: od: od: od:
   L: end:
seq(numtheory:-sigma(n), n = 1 .. 25);
seq(nops(TestSigma(n)), n = 1 .. 25);
\end{Verbatim}

We can experimentally verify Theorem \ref{ThmSigma} by extending the base code with the following one.

\begin{Verbatim}
k := 3:
t := n -> n+1:
w := n -> n+7:
M := n ->
   C(n^2, k, t(n), w(n))
   + A(1, [2, 0, 2], [2, 2, 2], [0, 0, 0], k, t(n), w(n))
   + A(1, [0, 2, 2], [2, 2, 2], [0, 0, 0], k, t(n), w(n))
   + A(1, [0, 0, 2], [2, 2, 2], [0, 0, 0], k, t(n), w(n))
   + A(-2*n, [1, 0, 1], [2, 2, 2], [0, 0, 0], k, t(n), w(n))
   + A(-2*n, [0, 1, 1], [2, 2, 2], [0, 0, 0], k, t(n), w(n))
   + A(-2*n, [0, 0, 1], [2, 2, 2], [0, 0, 0], k, t(n), w(n))
   + A(2, [1, 1, 2], [2, 2, 2], [0, 0, 0], k, t(n), w(n))
   + A(2, [1, 0, 2], [2, 2, 2], [0, 0, 0], k, t(n), w(n))
   + A(2, [0, 1, 2], [2, 2, 2], [0, 0, 0], k, t(n), w(n)):
Sigma := n -> HW(M(n))/w(n)-t(n)^k:
seq(numtheory:-sigma(n), n = 1 .. 25);
seq(Sigma(n), n = 1 .. 25);
\end{Verbatim}

We can use the following code in order to test Lemma \ref{LemmaSetPhi}.

\begin{Verbatim}
TestPhi := proc(n)
   local a, b, c, L:
   L := []:
   for a from 0 to n do
      for b from 0 to n do
         for c from 0 to n do
            if a*b-c*n-1 = 0
            then L := [op(L), [a, b, c]]: fi: od: od: od:
   L: end:
seq(numtheory:-phi(n), n = 1 .. 25);
seq(nops(TestPhi(n)), n = 1 .. 25);
\end{Verbatim}

We can experimentally verify Theorem \ref{ThmPhi} by extending the base code with the following one.

\begin{Verbatim}
k := 3:
t := n -> n+1:
w := n -> n+5:
M := n ->
   C(1, k, t(n), w(n))
   + A(-2*n, [1, 1, 1], [2, 2, 2], [0, 0, 0], k, t(n), w(n))
   + A(1, [2, 2, 0], [2, 2, 2], [0, 0, 0], k, t(n), w(n))
   + A(-2, [1, 1, 0], [2, 2, 2], [0, 0, 0], k, t(n), w(n))
   + A(n^2, [0, 0, 2], [2, 2, 2], [0, 0, 0], k, t(n), w(n))
   + A(2*n, [0, 0, 1], [2, 2, 2], [0, 0, 0], k, t(n), w(n)):
Phi := n -> HW(M(n))/w(n)-t(n)^k:
seq(numtheory:-phi(n), n = 1 .. 25);
seq(Phi(n), n = 1 .. 25);
\end{Verbatim}

We can use the following code in order to test Lemma \ref{LemmaSetInv}.

\begin{Verbatim}
TestInv := proc(m, n)
   local a, b, c, d, L:
   L := []:
   for a from 0 to n do
      for b from 0 to n do
         for c from 0 to n do
            for d from 0 to n do
               if (m*a-n*b-1)^2+(a-c-d-1)^2 = 0
               then L := [op(L), [a, b, c, d]]: fi: od: od: od: od:
   L: end:
Coprimes := n -> {select(i -> (gcd(i, n) = 1), [$2 .. n-1])[]}:
for n from 2 to 15 do seq(1/m mod n, m in Coprimes(n)): od;
for n from 2 to 15 do seq(nops(TestInv(m, n)), m in Coprimes(n)): od;
\end{Verbatim}

We can experimentally verify Theorem \ref{ThmInv} by extending the base code with the following one.

\begin{Verbatim}
k := 4:
t := n -> n+1:
w := n -> n+5:
M := (m, n) ->
   C(2, k, t(n), w(n))
   + A(m^2+1, [2, 0, 0, 0], [2, 2, 2, 2], [0, 0, 0, 0], k, t(n), w(n))
   + A(-2*(m+1), [1, 0, 0, 0], [2, 2, 2, 2], [0, 0, 0, 0], k, t(n), w(n))
   + A(n^2, [0, 2, 0, 0], [2, 2, 2, 2], [0, 0, 0, 0], k, t(n), w(n))
   + A(2*n, [0, 1, 0, 0], [2, 2, 2, 2], [0, 0, 0, 0], k, t(n), w(n))
   + A(-2*m*n, [1, 1, 0, 0], [2, 2, 2, 2], [0, 0, 0, 0], k, t(n), w(n))
   + A(2, [0, 0, 1, 1], [2, 2, 2, 2], [0, 0, 0, 0], k, t(n), w(n))
   + A(2, [0, 0, 1, 0], [2, 2, 2, 2], [0, 0, 0, 0], k, t(n), w(n))
   + A(2, [0, 0, 0, 1], [2, 2, 2, 2], [0, 0, 0, 0], k, t(n), w(n))
   + A(-2, [1, 0, 1, 0], [2, 2, 2, 2], [0, 0, 0, 0], k, t(n), w(n))
   + A(-2, [1, 0, 0, 1], [2, 2, 2, 2], [0, 0, 0, 0], k, t(n), w(n))
   + A(1, [0, 0, 2, 0], [2, 2, 2, 2], [0, 0, 0, 0], k, t(n), w(n))
   + A(1, [0, 0, 0, 2], [2, 2, 2, 2], [0, 0, 0, 0], k, t(n), w(n)):
Inv := (m, n) -> HW(M(m, n))/w(n)-t(n)^k:
Coprimes := n -> {select(i -> (gcd(i, n) = 1), [$2 .. n-1])[]}:
for n from 2 to 10 do seq(1/m mod n, m in Coprimes(n)): od;
for n from 2 to 10 do seq(Inv(m, n), m in Coprimes(n)): od;
\end{Verbatim}

We can use the following code in order to test Lemma \ref{LemmaSetRoot}.

\begin{Verbatim}
TestRoot := proc(m, n)
   local a, b, L:
   L := []:
   for a from 0 to n do
      for b from 0 to n do
         if a+b^m-n = 0
         then L := [op(L), [a, b]]: fi: od: od:
   L: end:
for m from 2 to 5 do seq(floor(n^(1/m))+1, n = 1 .. 35): od;
for m from 2 to 5 do seq(nops(TestRoot(m, n)), n = 1 .. 35): od;
\end{Verbatim}

We can experimentally verify Theorem \ref{ThmSquareRoot} by extending the base code with the following one.

\begin{Verbatim}
k := 4:
t := n -> n+1:
w := n -> n+5:
M := n ->
   C(n^2, k, t(n), w(n))
   + A(-2*n, [1, 0, 0, 0], [2, 2, 2, 2], [0, 0, 0, 0], k, t(n), w(n))
   + A(-2*n, [0, 0, 0, 1], [2, 2, 2, 2], [0, 0, 0, 0], k, t(n), w(n))
   + A(1, [0, 2, 2, 0], [2, 2, 2, 2], [0, 0, 0, 0], k, t(n), w(n))
   + A(1, [2, 0, 0, 0], [2, 2, 2, 2], [0, 0, 0, 0], k, t(n), w(n))
   + A(1, [0, 2, 0, 0], [2, 2, 2, 2], [0, 0, 0, 0], k, t(n), w(n))
   + A(1, [0, 0, 2, 0], [2, 2, 2, 2], [0, 0, 0, 0], k, t(n), w(n))
   + A(2, [0, 0, 0, 2], [2, 2, 2, 2], [0, 0, 0, 0], k, t(n), w(n))
   + A(2, [1, 0, 0, 1], [2, 2, 2, 2], [0, 0, 0, 0], k, t(n), w(n))
   + A(-2, [0, 1, 1, 1], [2, 2, 2, 2], [0, 0, 0, 0], k, t(n), w(n))
   + A(-2, [0, 1, 1, 0], [2, 2, 2, 2], [0, 0, 0, 0], k, t(n), w(n)):
Sqrt := n -> HW(M(n))/w(n)-t(n)^k-1:
seq(floor(sqrt(n)), n = 1 .. 15);
seq(Sqrt(n), n = 1 .. 15);
\end{Verbatim}

We can use the following code in order to test Lemma \ref{LemmaSetLog}.

\begin{Verbatim}
TestLog := proc(m, n)
   local a, b, L:
   L := []:
   for a from 0 to n do
      for b from 0 to n do
         if a+m^b-n = 0
         then L := [op(L), [a, b]]: fi: od: od:
   L: end:
for m from 2 to 5 do seq(floor(log[m](n))+1, n = 1 .. 25): od;
for m from 2 to 5 do seq(nops(TestLog(m, n)), n = 1 .. 25): od;
\end{Verbatim}

We can experimentally verify Theorem \ref{ThmLog} by extending the base code with the following one.

\begin{Verbatim}
k := 2:
t := n -> n+1:
w := (m, n) -> 2*m*n:
M := (m, n) ->
   C(n^2, k, t(n), w(m, n))
   + A(-2*n, [0, 0], [m, m], [0, 1], k, t(n), w(m, n))
   + A(-2*n, [1, 0], [m, m], [0, 0], k, t(n), w(m, n))
   + A(2, [1, 0], [m, m], [0, 1], k, t(n), w(m, n))
   + A(1, [0, 0], [m, m], [0, 2], k, t(n), w(m, n))
   + A(1, [2, 0], [m, m], [0, 0], k, t(n), w(m, n)):
Log := (m, n) -> HW(M(m, n))/w(m, n)-t(n)^k-1:
for m from 2 to 5 do seq(floor(log[m](n)), n = 1 .. 25): od;
for m from 2 to 5 do seq(Log(m, n), n = 1 .. 25): od;
\end{Verbatim}

We can use the following code in order to test Lemma \ref{LemmaSetOrd} (we skipped the case $ ( m , n ) = ( 4 , 5 ) $ because its computation takes too much time).

\begin{Verbatim}
with(numtheory):
TestOrd := proc(m, n)
   local a, b, c, d, L:
   L := []:
   for a from 0 to m^phi(n) do
      for b from 0 to m^phi(n) do
         for c from 0 to m^phi(n) do
            for d from 0 to m^phi(n) do
               if (m^a-n*b-1)^2+(a-c-1)^2+(phi(n)-d-a)^2 = 0
               then L := [op(L), [a, b, c, d]]: fi: od: od: od: od:
   L: end:
for n from 3 to 6 do seq(phi(n)/order(m, n),
m in `minus`({select(i -> (gcd(i, n) = 1), [$2 .. n-1])[]}, {4})): od;
for n from 3 to 6 do seq(nops(TestOrd(m, n)),
m in `minus`({select(i -> (gcd(i, n) = 1), [$2 .. n-1])[]}, {4})): od;
\end{Verbatim}

We can experimentally verify Theorem \ref{ThmOrd} by extending the base code with the following one (as the computations are very heavy at this point, we only check the identity $ \ord ( 2 , 5 ) = 4 $).

\begin{Verbatim}
with(numtheory):
m, n, k := 2, 5, 4:
t := m^phi(n)+1:
w := 2*m^(phi(n)+1)+2:
M :=
   C((phi(n))^2+2, k, t, w)
   + A(-2*n, [0, 1, 0, 0], [m, m, m, m], [1, 0, 0, 0], k, t, w)
   + A(-2, [0, 0, 0, 0], [m, m, m, m], [1, 0, 0, 0], k, t, w)
   + A(1, [0, 0, 0, 0], [m, m, m, m], [2, 0, 0, 0], k, t, w)
   + A(2*n, [0, 1, 0, 0], [m, m, m, m], [0, 0, 0, 0], k, t, w)
   + A(2, [0, 0, 1, 0], [m, m, m, m], [0, 0, 0, 0], k, t, w)
   + A(-2, [1, 0, 1, 0], [m, m, m, m], [0, 0, 0, 0], k, t, w)
   + A(2, [1, 0, 0, 1], [m, m, m, m], [0, 0, 0, 0], k, t, w)
   + A(2, [2, 0, 0, 0], [m, m, m, m], [0, 0, 0, 0], k, t, w)
   + A(n^2, [0, 2, 0, 0], [m, m, m, m], [0, 0, 0, 0], k, t, w)
   + A(1, [0, 0, 2, 0], [m, m, m, m], [0, 0, 0, 0], k, t, w)
   + A(1, [0, 0, 0, 2], [m, m, m, m], [0, 0, 0, 0], k, t, w)
   + A(-2*(phi(n)+1), [1, 0, 0, 0], [m, m, m, m], [0, 0, 0, 0], k, t, w)
   + A(-2*phi(n), [0, 0, 0, 1], [m, m, m, m], [0, 0, 0, 0], k, t, w):
phi(n)/(HW(M)/w-t^k);
\end{Verbatim}

We can use the following code in order to test Lemma \ref{LemmaSetDiscLog} (for the case $ n = 5 $, which already takes some time).

\begin{Verbatim}
with(numtheory):
TestDiscLog := proc(m, g, n)
    local a, b, c, d, L:
    L := []:
    for a from 0 to g^phi(n) do
        for b from 0 to g^phi(n) do
            for c from 0 to g^phi(n) do
                for d from 0 to g^phi(n) do
                    if (a+b+c+1-phi(n))^2+(g^(a+b+1)-n*d-m)^2 = 0
                    then L := [op(L), [a, b, c, d]]: fi: od: od: od: od:
    L: end:
n := 5:
Coprimes := {select(i -> (gcd(i, n) = 1), [$2 .. n-1])[]}:
PrimitiveRoots := {select(i -> (i in Coprimes and
    order(i, n) = phi(n)), [$2 .. n-1])[]}:
for g in PrimitiveRoots do seq(mlog(m, g, n), m in Coprimes): od;
for g in PrimitiveRoots do seq(nops(TestDiscLog(m, g, n)), m in Coprimes): od;
\end{Verbatim}

We can experimentally verify Theorem \ref{ThmDiscLog} (for the case $ n = 4 $) by extending the base code with the following one.

\begin{Verbatim}
with(numtheory):
k := 4:
t := (g, n) -> g^phi(n)+1:
w := (g, n) -> 5+2*n+2*g*(2*g^phi(n)+1):
M := (m, g, n) ->
    C(1+m^2-2*phi(n)+(phi(n))^2, k, t(g, n), w(g, n))
    + A(2*(1-phi(n)), [1, 0, 0, 0], [g, g, g, g],
        [0, 0, 0, 0], k, t(g, n), w(g, n))
    + A(2*(1-phi(n)), [0, 1, 0, 0], [g, g, g, g],
        [0, 0, 0, 0], k, t(g, n), w(g, n))
    + A(2*(1-phi(n)), [0, 0, 1, 0], [g, g, g, g],
        [0, 0, 0, 0], k, t(g, n), w(g, n))
    + A(-2*m*g, [0, 0, 0, 0], [g, g, g, g],
        [1, 1, 0, 0], k, t(g, n), w(g, n))
    + A(g^2, [0, 0, 0, 0], [g, g, g, g],
        [2, 2, 0, 0], k, t(g, n), w(g, n))
    + A(-2*g*n, [0, 0, 0, 1], [g, g, g, g],
        [1, 1, 0, 0], k, t(g, n), w(g, n))
    + A(1, [2, 0, 0, 0], [g, g, g, g],
        [0, 0, 0, 0], k, t(g, n), w(g, n))
    + A(1, [0, 2, 0, 0], [g, g, g, g],
        [0, 0, 0, 0], k, t(g, n), w(g, n))
    + A(1, [0, 0, 2, 0], [g, g, g, g],
        [0, 0, 0, 0], k, t(g, n), w(g, n))
    + A(2*m*n, [0, 0, 0, 1], [g, g, g, g],
        [0, 0, 0, 0], k, t(g, n), w(g, n))
    + A(n^2, [0, 0, 0, 2], [g, g, g, g],
        [0, 0, 0, 0], k, t(g, n), w(g, n))
    + A(2, [1, 1, 0, 0], [g, g, g, g],
        [0, 0, 0, 0], k, t(g, n), w(g, n))
    + A(2, [1, 0, 1, 0], [g, g, g, g],
        [0, 0, 0, 0], k, t(g, n), w(g, n))
    + A(2, [0, 1, 1, 0], [g, g, g, g],
        [0, 0, 0, 0], k, t(g, n), w(g, n)):
DiscLog := (m, g, n) -> HW(M(m, g, n))/w(g, n)-t(g, n)^k:
mlog(3, 3, 4);
DiscLog(3, 3, 4);
\end{Verbatim}

\end{appendices}


\begin{thebibliography}{99}

\bibitem{BorweinCrandall}
\newblock J.\ M.\ Borwein and R.\ E.\ Crandall,
\newblock Closed Forms: What They Are and Why We Care,
\newblock \textit{Notices of the American Mathematical Society} \textbf{60} (2013).
\newblock \url{https://doi.org/10.1090/NOTI936}.

\bibitem{CrandallPomerance}
\newblock R.\ Crandall and C.\ Pomerance,
\newblock \textit{Prime Numbers: A Computational Perspective} (2nd ed.),
\newblock Springer,
\newblock 2005.

\bibitem{Enderton}
\newblock H.\ B.\ Enderton,
\newblock \textit{A Mathematical Introduction to Logic} (2nd ed.),
\newblock Academic Press,
\newblock 2001.

\bibitem{HardyWright}
\newblock G.\ H.\ Hardy and E.\ M.\ Wright,
\newblock \textit{An Introduction to the Theory of Numbers} (6th ed.),
\newblock Oxford University Press (2008).

\bibitem{KrizekEtAl}
\newblock M.\ K\v{r}\'{i}\v{z}ek, L.\ Somer and A.\ \v{S}olková,
\newblock \textit{From Great Discoveries in Number Theory to Applications},
\newblock Springer Cham,
\newblock 2021.

\bibitem{Marchenkov}
\newblock S.\ S.\ Marchenkov,
\newblock Superpositions of Elementary Arithmetic Functions,
\newblock \textit{Journal of Applied and Industrial Mathematics} \textbf{3} (2007).

\bibitem{Matiyasevich}
\newblock Y.\ Matiyasevich,
\newblock \textit{Hilbert's Tenth Problem},
\newblock MIT Press, Sprin\-ger-Ver\-lag, New York,
\newblock 1993.

\bibitem{Mazzanti}
\newblock S.\ Mazzanti,
\newblock Plain Bases for Classes of Primitive Recursive Functions,
\newblock \textit{Mathematical Logic Quarterly} \textbf{48} (2002).

\bibitem{Mendelson}
\newblock E.\ Mendelson,
\newblock \textit{Introduction to Mathematical Logic} (6th ed.),
\newblock Taylor \& Francis,
\newblock 2015.

\bibitem{Nederlof}
\newblock J.\ Nederlof,
\newblock \textit{Space and Time Efficient Structural Improvements of Dynamic Programming Algorithms},
\newblock doctoral thesis,
\newblock University of Bergen,
\newblock 2011.

\bibitem{Oitavem}
\newblock I.\ Oitavem,
\newblock New recursive characterizations of the elementary functions and the functions computable in polynomial space,
\newblock \textit{Revista Matemática de la Universidad Complutense de Madrid} \textbf{10} (1997).

\bibitem{PetkovsekEtAl}
\newblock M.\ Petkovšek, H.\ S.\ Wilf, D.\ Zeilberger,
\newblock \textit{A = B},
\newblock A.\ K.\ Peters / CRC Press,
\newblock 1996.

\bibitem{PrunescuSaurasAltuzarra}
\newblock M.\ Prunescu and L.\ Sauras-Altuzarra,
\newblock An arithmetic term for the factorial function,
\newblock \textit{Examples \& Counterexamples} \textbf{5} (2024).
\newblock \url{https://doi.org/10.1016/j.exco.2024.100136}.

\bibitem{Rosen}
\newblock K.\ H.\ Rosen,
\newblock \textit{Elementary Number Theory and Its Applications} (6th ed.),
\newblock Addison-Wesley,
\newblock 2011.

\bibitem{SaurasAltuzarra}
\newblock L.\ Sauras-Altuzarra,
\newblock \textit{From Logic to Discrete Geometry via Lattices},
\newblock doctoral thesis,
\newblock Vienna University of Technology,
\newblock 2024.
\newblock \url{https://doi.org/10.34726/hss.2024.111390}

\bibitem{SaurasAltuzarra2}
\newblock L.\ Sauras-Altuzarra,
\newblock \textit{Hypergeometric closed forms},
\newblock master thesis,
\newblock Vienna University of Technology,
\newblock 2018.
\newblock \url{https://doi.org/10.25365/thesis.57260}

\bibitem{Stein}
\newblock W.\ Stein,
\newblock \textit{Elementary Number Theory: Primes, Congruences, and Secrets},
\newblock Springer-Verlag,
\newblock 2017.

\bibitem{Encyclopedia}
\newblock Various,
\newblock ``P-adic valuation'',
\newblock Encyclopedia of Mathematics.
\newblock \url{https://encyclopediaofmath.org/index.php?title=P-adic_valuation}

\bibitem{VereschchaginShen}
\newblock N.\ K.\ Vereschchagin and A.\ Shen,
\newblock \textit{Computable Functions} (translated by V.\ N.\ Dubrovskii),
\newblock American Mathematical Society,
\newblock 2002.

\bibitem{Weisstein0}
\newblock E.\ W.\ Weisstein,
\newblock ``Pairing Function'',
\newblock from MathWorld -- A Wolfram Web Resource.
\newblock \url{https://mathworld.wolfram.com/PairingFunction.html}

\bibitem{Weisstein1}
\newblock E.\ W.\ Weisstein,
\newblock ``RSA Encryption'',
\newblock from MathWorld -- A Wolfram Web Resource.
\newblock \url{https://mathworld.wolfram.com/RSAEncryption.html}

\bibitem{Weisstein2}
\newblock E.\ W.\ Weisstein,
\newblock ``Semiprime'',
\newblock from MathWorld -- A Wolfram Web Resource.
\newblock \url{https://mathworld.wolfram.com/Semiprime.html}

\end{thebibliography}
\end{document}